\documentclass[a4paper,12pt]{article}
\usepackage[latin1]{inputenc}
\usepackage[T1]{fontenc}
\usepackage{amsmath,amssymb,bm}
\usepackage[2001]{abreviations}
\usepackage[2001,numero]{theo}
\usepackage{fleches}
\usepackage[english,frenchb]{babel}
\def\english{\shorthandoff{;:!?}}
\usepackage{url}
\usepackage[all,dvips,poly]{xy}



\usepackage{changebar}
\let\cbstart\relax
\let\cbend\relax

\let\wh\widehat
\let\wt\widetilde
\def\et{_{\text{\rm{ét}}}}
\let\d\textbf
\DeclareMathOperator{\pr}{pr}
\DeclareMathOperator{\desc}{\cat{Desc}}
\DeclareMathOperator{\LC}{LC}
\DeclareMathOperator{\SLC}{SLC}
\DeclareMathOperator{\SLCF}{SLCF}
\def\T{{\mathscr T}}
\def\B{{\cat B}}
\usepackage{ifthen}
\newcommand{\F}[1]{\relax
  \ifthenelse{\equal{#1}{1}}%
  {\cat{F}^{\mathstrut}_X}{%
    \ifthenelse{\equal{#1}{2}}%
    {\cat{F}^{\mathstrut}_{X \Times_Y X}}%
    {\cat{F}^{\mathstrut}_{X \Times_Y X \Times_Y X}}%
    }}      

\qsdfleche{tlra}

\makeatletter
  \def\dsa{\mathbin{\overset{2}{*}}}
  \def\Dsa{\@ifnextchar_\Dsa@i\dsa}
  \def\Dsa@i_#1{\mathbin{\mathop{\overset{2}{*}}\limits_{#1}}}
  \def\TIMES#1#2#3{\mathbin{\mathop{
        {}_{#1}^{\hphantom{#3}}{\times}{}_{#3}^{\hphantom{#1}}
        }\limits_{#2}}}
  \def\dtimes{\mathbin{\mathop{\times}\limits^2}}
  \def\dtimes{\mathchoice
    {\mathbin{\smash{\mathop{\times}\limits^{
            \hbox{\smash{\raisebox{-.7ex}{\ensuremath{\scriptscriptstyle 2}}}}
            }}}}%
    {\mathbin{\smash{\mathop{\times}\limits^{
            \hbox{\smash{\raisebox{-.7ex}{\ensuremath{\scriptscriptstyle 2}}}}
            }}}}%
    {\mathbin{\mathop{\times}\limits^{2}}}%
    {\mathbin{\mathop{\times}\limits^{2}}}%
    }
  \let\OLD@dtimes\dtimes
  \def\DTimes{\@ifnextchar_\dtimesWSB\OLD@dtimes}
  \def\dtimesWSB_#1{\mathbin{\mathop{\OLD@dtimes}\limits_{#1}}}
\makeatother

\newcommand{\UUN}[4][r]{%
  \ar@/^2pc/[#1]^{#2}_*=<0.3pt>{}="HAUT"
  \ar@/_2pc/[#1]_{#3}^*=<0.3pt>{}="BAS"
  \save\POS "HAUT",*{
      \vrule height 2pt depth 2pt width 0pt
      \vrule height 0pt depth 0pt width 4pt
      }="HAUT",\restore
  \save\POS "BAS",*{
      \vrule height 2pt depth 2pt width 0pt
      \vrule height 0pt depth 0pt width 4pt
      }="BAS",\restore
  \ar @{=>} "HAUT";"BAS" ^{#4}
}

\DeclareMathOperator{\LCF}{LCF}
\def\cartesien{%
  \ar@{-}[]+R+<6pt,-2pt>;[]+RD+<6pt,-6pt>%
  \ar@{-}[]+D+<2pt,-6pt>;[]+RD+<6pt,-6pt>%
}
\def\cocartesien{%
  \ar@{-}[]+L+<-6pt,+2pt>;[]+LU+<-6pt,+6pt>%
  \ar@{-}[]+U+<-2pt,+6pt>;[]+LU+<-6pt,+6pt>%
}
\def\enlargexyentry#1{%
  \POS "#1",*{
    \vrule height 2pt depth 2pt width 0pt
    \vrule height 0pt depth 0pt width 4pt
    }="#1",}

\makeatletter
\def\bighole{\hbox{\dimen@=3\objectmargin@ \kern2\dimen@
  \vrule height\dimen@ depth\dimen@ width\z@}}
\makeatother

\newsavebox\tempboxa
\newsavebox\tempboxb
\newsavebox\tempboxc
\newsavebox\tempboxd
\newsavebox\tempboxe

\title{\large THÉORÈME DE VAN KAMPEN \\ POUR LES CHAMPS ALGÉBRIQUES}
\author{\normalsize Vincent \textsc{Zoonekynd} 
  \\ \footnotesize \url{<zoonek@math.jussieu.fr>}}
\date{}
\begin{document}
\maketitle

\selectlanguage{english}
\begin{abstract}
\cbstart
  We define a category whose objects are finite
  étale coverings of an algebraic stack and prove
  that it is a Galois category and that it allows one to compute the
  fundamental group of the stack. 

  We then prove a Van Kampen theorem for algebraic
  stacks whose simplest form reads:
  Let $U$ and $V$ be open substacks of an
  algebraic stack $X$ with $X = U \union V$,\ifmmode\,\fi{}  let
  $P$ be a set of base points,\ifmmode\,\fi{}  at least one in
  each connected component of $X$,\ifmmode\,\fi{}  $U$,\ifmmode\,\fi{}  $V$ and $U
  \inter V$,\ifmmode\,\fi{}  then there is a cocartesian square of
  fundamental progroupoids
  $$ \xymatrix{
    \pi_1( U \inter V ,\ifmmode\,\fi{}  P) \ar[r] \ar[d] & 
    \pi_1(U,\ifmmode\,\fi{} P) \ar[d] \\
    \pi_1(V,\ifmmode\,\fi{} P) \ar[r] &
    \pi_1( X ,\ifmmode\,\fi{}  P). \cocartesien }$$
\cbend
\end{abstract}
\selectlanguage{frenchb}
  

\section*{Introduction}
\cbstart
Le groupe fondamental d'un espace topologique
(semi-localement simplement connexe) peut se
définir de deux manières : à l'aide de chemins
tracés sur cet espace ou à l'aide de ses
revêtements. Cette deuxième définition se
généralise aisément aux schémas : la catégorie des
revêtements finis d'un schéma $X$,\ifmmode\,\fi{}  \ie,\ifmmode\,\fi{}  la
catégorie dont les objets sont les morphismes
étales finis $Y \lra X$ et dont les morphismes 
$(Y\lra X) \lra (Y' \lra X)$ sont les diagrammes
commutatifs 
$$ \xymatrix{ Y \ar[rr] \ar[rd] && Y' \ar[ld] \\ &
  X } $$
est galoisienne,\ifmmode\,\fi{}  \ie,\ifmmode\,\fi{}  équivalente à la catégorie
des ensembles finis munis d'une action continue
d'un certain groupe profini,\ifmmode\,\fi{}  uniquement déterminé (à
isomorphisme non unique près) ; c'est ce groupe
que l'on définit comme le groupe fondamental
profini du schéma $X$. 

Nous nous intéressons à la généralisation de ces
résultats aux champs algébriques \cite{LMB,vistoli}, qui
apparaissent par exemple dans l'étude des espaces
de modules de courbes : on peut les voir comme des
« objets géométriques » qui « ressemblent
localement »,\ifmmode\,\fi{}  pour la topologie étale,\ifmmode\,\fi{}  à des
quotients de schémas par des groupes finis,\ifmmode\,\fi{}  ou
comme des analogues algébriques des orbifolds 
\cite{thurston,MP}.
On peut définir le groupe fondamental d'un champ
algébrique à l'aide de schémas simpliciaux 
\cite{oda,friedlander} où à l'aide de faisceaux
étales localement constants finis
\cite{leroy,moerdijk,Z:pi1} qui jouent le rôle des
revêtements. Une définition du groupe fondamental
d'un champ algébrique à l'aide d'une catégorie de
revêtements est problématique,\ifmmode\,\fi{}  car les champs
algébriques ne forment pas tant une catégorie
qu'une 2-catégorie. Nous montrons comment
définir cette 1-catégorie de revêtements 
(définitions 
\ref{defi:dcataudessus},\ifmmode\,\fi{} 
\ref{def:revchamp})
et
montrons qu'elle est galoisienne et équivalente à
la catégorie des faisceaux étales localement
constants finis 
(théorème \ref{theo:revchamp})
: elle permet donc de calculer le groupe fondamental d'un champ. 

Nous nous intéressons ensuite à l'un des moyens de
calculer des groupes fondamentaux : le théorème de
Van Kampen. Sous sa forme classique,\ifmmode\,\fi{}  il permet de
calculer le groupe fondamental d'une réunion $U
\union V$,\ifmmode\,\fi{}  connaissant celui des morceaux $U$,\ifmmode\,\fi{}  $V$
et de leur intersection $U \inter V$,\ifmmode\,\fi{}  sous réserve
tout soit connexe. Il peut aussi se formuler en
termes de \emph{groupoïdes} \cite{brown}, ce qui
permet de se libérer des hypothèses de connexité.
Le théorème de Van Kampen est déjà connu pour les
orbifolds : il est par exemple démontré par
\cite{H} à l'aide d'espaces classifiants ou par
\cite[4.13]{pronk} à l'aide de la propriété universelle
définissant le groupoïde fondamental. Quant-à
nous,\ifmmode\,\fi{}  c'est à l'aide de la théorie de la descente
\cite{giraud} que nous établirons un théorème de
Van Kampen pour le groupoïde fondamental d'un champs algébrique,\ifmmode\,\fi{}  et
retrouverons donc le résultat pour les orbifolds
(théorème \ref{theo:VKchamps}).

Voici maintenant le plan de ce travail.

Dans une première section,\ifmmode\,\fi{}  nous rappelons la
définition du groupe fondamental d'un champ
algébrique à l'aide de faisceaux étales localement
constants finis \cite{leroy,moerdijk,Z:pi1,zoonek}.

Dans une seconde section,\ifmmode\,\fi{}  nous considérons une
2-catégorie $\dcat C$ et un objet 
$X \in \ob \dcat C$ et nous définissons la
2-catégorie $\dcat C / X$ des objets de $\dcat C$
au dessus de~$X$,\ifmmode\,\fi{}  puis la 1-catégorie associée ;
nous montrons que ces notions se comportent « bien
» vis-à-vis des produits fibrés. 

Dans une troisième section,\ifmmode\,\fi{}  nous utilisons ces
constructions pour définir la 2-catégorie des
revêtements étales finis d'un champ algébrique
$X$,\ifmmode\,\fi{}  comme une sous-catégorie pleine de 
la 2-catégorie $\dcat{Champs}/X$ des champs au dessus de
$X$,\ifmmode\,\fi{}  puis montrons que la 1-catégorie associée est
équivalente à la catégorie des faisceaux étales
localement constants finis sur $X$ : elle permet
de calculer le groupe fondamental. 

Dans une quatrième section,\ifmmode\,\fi{}  nous énonçons un
théorème de Van Kampen pour le calcul du groupoïde
fondamental d'un topos dans lequel on a choisi des
points-base et le démontrons à l'aide du théorème
de Van Kampen « usuel » exprimé en termes de
données de descente.

Dans une cinquième et dernière section,\ifmmode\,\fi{}  nous
appliquons ce théorème au cas des champs
algébriques. 
\cbend

\section{Groupoïde fondamental d'un champ\\algébrique}
Nous rappelons brièvement les résultats de
\cite{leroy} et \cite{zoonek,Z:pi1} définissant le
progroupoïde fondamental d'un topos localement
connexe,\ifmmode\,\fi{}  en particulier du topos des faisceaux
étales sur un champ algébrique,\ifmmode\,\fi{}  à partir de la
sous-catégorie des objets localement constants,\ifmmode\,\fi{} 
qui est équivalente au topos classifiant d'un
progroupoïde.

\begin{notations}
  \cbstart 
  Si $\T$ est un topos,\ifmmode\,\fi{}  nous noterons
  $\nothing$ son objet initial et~$*$ son objet
  final.  
  \cbend
\end{notations}

\begin{defi}
  Un objet $X$ d'un topos $\T$ est \d{connexe} si
  pour tout isomorphisme $X \iso A \amalg B$,\ifmmode\,\fi{}  on a
  $A \iso \nothing$ ou $B \iso \nothing$.
\end{defi}

\begin{defi}
  Un topos est \d{localement connexe} s'il est
  engendré par ses objets connexes ou,\ifmmode\,\fi{}  ce qui
  revient au même,\ifmmode\,\fi{}  si le foncteur \d{objet
    constant}
  $$
  \definefunction{\cat{Ens}}{\T}{I}{I_\T =
    \coprod_{i \in I} *} $$
  (où $\cat{Ens}$
  désigne la catégorie des ensembles) possède un
  adjoint à gauche,\ifmmode\,\fi{}  noté $\pi$ et appelé foncteur
  « composantes connexes ».
\end{defi}

\begin{defi}
  Un objet $F$ d'un topos $\T$ est \d{trivialisé}
  par un objet $U$ si sa restriction $F |_U = ( F
  \times U \lra U)$ est un objet constant de
  $\T/U$.
\end{defi}

\begin{defi}
  Un objet $F$ d'un topos $\T$ est \d{localement
    constant} s'il existe un crible $R$ couvrant
  l'objet final $*$ tel que tous les objets
  connexes de~$R$ trivialisent $F$.
\end{defi}

\begin{defi}
  \cbstart
  Un \d{point} d'un topos $\T$ est un morphisme de
  topos depuis le topos ponctuel (la catérorie des
  ensembles),\ifmmode\,\fi{}  \ie,\ifmmode\,\fi{}  
  $\cat{Ens} \lra \T$.

  Les points d'un topos $\T$ et leurs
  isomorphismes forment un groupoïde que nous
  noterons 
  $\cat{Points}\,\ifmmode\,\fi{} \T$.
  
  Si $G$ est un groupoïde,\ifmmode\,\fi{}  nous noterons $\cat BG$
  son \d{topos classifiant},\ifmmode\,\fi{}  \ie,\ifmmode\,\fi{}  la catégorie des
  préfaisceaux sur le groupoïde opposé à $G$. Par
  exemple,\ifmmode\,\fi{}  si $G$ est un groupe,\ifmmode\,\fi{}  $\cat BG$ est la
  catégorie des ensembles munis d'une action de
  $G$ \emph{à gauche}.  
  \cbend
\end{defi}
\begin{theo}\label{theo:pi1topos}%
  Soit $\T$ un topos localement connexe. 
  
  Pour tout crible $R$ couvrant l'objet final $*$
  de $\T$,\ifmmode\,\fi{}  la sous-catégorie pleine $\LC(\T\!,\ifmmode\,\fi{} R)$
  des objets localement constants trivialisés par
  les objets connexes de $R$ est équivalente au
  topos classifiant du groupoïde
  $\cat{Points}\,\ifmmode\,\fi{} \LC(\T\!,\ifmmode\,\fi{} R)$ des points de
  $\LC(\T\!,\ifmmode\,\fi{} R)$.
  
  La sous-catégorie pleine $\SLC \T$ des somme
  disjointes d'objets localement constants de $\T$
  est un topos qui s'identifie à la $2$-limite
  projective des topos classifiants de groupoïdes
  $\LC(\T\!,\ifmmode\,\fi{} R)$,\ifmmode\,\fi{} 
  $$ \SLC\T = \dlimpro_{R \in J(*)} \LC(\T\!,\ifmmode\,\fi{} R). $$
\end{theo}

\begin{dem}
  Voir \cite{leroy}.
\end{dem}

\begin{defi}
  Un \d{point-base} d'un topos localement connexe $\T$ est un
  morphisme de topos 
  $\cat{Ens} \lra \SLC \T$. 
\end{defi}

\begin{rem}
  Un point-base de $\T$ définit une composante
  connexe de $\SLC\T$ et même de $\LC(\T\!,\ifmmode\,\fi{} R)$ : on
  peut le démontrer comme suit.  D'après
  \cite[I.1.11.a]{zoonek}, les composantes
  connexes de $\SLC\T$ sont exactement celles de
  $\LC(\T\!,R)$.  Or, d'après \cite[3.2.8]{leroy},
  $$\LC(\T\!,\ifmmode\,\fi{} R) \equiv \cat B\,\ifmmode\,\fi{}  \cat{Points}\,\ifmmode\,\fi{} 
  \LC(\T\!,\ifmmode\,\fi{} R),\ifmmode\,\fi{} $$
  où $\cat B \,\ifmmode\,\fi{}  \cat{Points}\,\ifmmode\,\fi{} 
  \LC(\T\!,\ifmmode\,\fi{} R)$ désigne le topos classifiant du
  groupoïde des points de $\LC(\T\!,\ifmmode\,\fi{}  R)$.  Donc un
  point de $\SLC\T$ définit un point de
  $\LC(\T\!,\ifmmode\,\fi{} R)$
  $$
  \cat{Ens} \lra \SLC\T \lra \LC(\T\!,\ifmmode\,\fi{} R),\ifmmode\,\fi{}  $$
  \ie,\ifmmode\,\fi{} 
  un objet du groupoïde $\cat{Points}\,\ifmmode\,\fi{} 
  \LC(\T\!,\ifmmode\,\fi{} R)$,\ifmmode\,\fi{}  donc une composante connexe de ce
  groupoïde,\ifmmode\,\fi{}  donc une composante connexe de son
  topos classifiant.

  Un choix de points-bases $P = (p_i :
  \cat{Ens} \lra \SLC\T)_{i\in I}$,\ifmmode\,\fi{}  avec au moins un point-base dans
  chaque composante de $\SLC\T$,\ifmmode\,\fi{}  permet de définir des groupoïdes 
  $\pi_1(\T\!,\ifmmode\,\fi{} R,\ifmmode\,\fi{} P)$,\ifmmode\,\fi{}  dont l'ensemble des objets est $P$,\ifmmode\,\fi{} 
  et des équivalences 
  $\LC(\T\!,\ifmmode\,\fi{} R) \equiv \cat B \,\ifmmode\,\fi{}  \pi_1(\T\!,\ifmmode\,\fi{} R,\ifmmode\,\fi{} P)$. 
  On a donc 
  $$ \cat B \,\ifmmode\,\fi{}  \proobj_{R \in J(*)} \pi_1(\T\!,\ifmmode\,\fi{} R,\ifmmode\,\fi{} P)  
  := \dlimpro_{R \in J(*)} \cat B \,\ifmmode\,\fi{}  \pi_1(\T\!,\ifmmode\,\fi{} R,\ifmmode\,\fi{} P)  
  \equiv \SLC\T. $$
  Nous dirons que 
  $\smash{\pi_1(\T\!,\ifmmode\,\fi{} P) = \dlimpro\limits_{R\in J(*)} \pi_1(\T\!,\ifmmode\,\fi{} R,\ifmmode\,\fi{} P)}$ 
  est le \d{progroupoïde fondamental} de $\T$. 
\end{rem}

\begin{notations}
  \cbstart
  Par \d{champ algébrique},\ifmmode\,\fi{}  nous entendrons
  toujours « champ de Deligne--Mumford » 
  \cite{DM,LMB,vistoli}. Les champs algébriques 
  forment une 2-catégorie 
  $\dcat{Champs}$ dont les 2-morphismes sont des
  isomorphismes. 
  \cbend
\end{notations}

\begin{defi}
  Le \d{site étale} $X\et$ d'un champ algébrique 
  $X$ a pour objets les morphismes 
  étales 
  $f : T \lra X$,\ifmmode\,\fi{}  où $T$ est un schéma,\ifmmode\,\fi{}  pour morphismes 
  $f \lra f'$ les couples 
  \cbstart
  $(\phi : T \lra T',\ifmmode\,\fi{} \alpha:f'\phi\lra f)$
  \cbend
  $$ \english 
  \xymatrix{ 
    T \ar[rr]^\phi \ar[rd]_f|*{}="A"
    \enlargexyentry A
    && T' \ar@{=>}"A"_-\alpha \ar[ld]^{f'} \\ & X
    }$$
  et pour familles couvrantes les familles épimorphiques.
\end{defi}

\begin{lemme}\label{lemme:deuxsites}%
  Soit $X$ un champ algébrique. 
  Considérons le site 
  $X'\et$,\ifmmode\,\fi{}  dont les objets sont les morphismes étales 
  $f : T \lra X$,\ifmmode\,\fi{}  où $T$ est un champ algébrique,\ifmmode\,\fi{}  dont les morphismes 
  sont les couples 
  \cbstart
  $(\phi : T \lra T',\ifmmode\,\fi{} \alpha:f'\phi\lra f)$
  \cbend
  $$ \english 
  \xymatrix{ 
    T \ar[rr]^\phi \ar[rd]_f|*{}="A"
    \enlargexyentry A
    && T' \ar@{=>}"A"_-\alpha \ar[ld]^{f'} \\ & X,
    }$$
  modulo la relation d'équivalence 
  \begin{align*}
    (\phi,\ifmmode\,\fi{} \alpha) \sim (\phi',\ifmmode\,\fi{} \alpha') &\Longleftrightarrow 
    \exists \beta 
    \qquad 
    \raisebox{.5\depth}{\xymatrix{
        T \ar[rd]_f \UN[rr]{\phi}{\phi'}{\beta} && 
        T' \ar[ld]^{f'} \\ & X 
        }}
    \text{ $2$-commutatif} \\
    &\Longleftrightarrow \exists \beta : \phi \lra \phi' 
    \text{ tel que }
    \raisebox{.5\depth}{\xymatrix{
        f'\phi \ar[rd]_\alpha \ar[rr]^{f'\beta} && f'\phi'
        \ar[ld]^{\alpha'} \\ & f }}
  \end{align*}
  et dont les familles couvrantes sont les familles épimorphiques.

  Les topos $\cat{Sh}\,\ifmmode\,\fi{} X\et$ et $\cat{Sh}\,\ifmmode\,\fi{} X'\et$ sont équivalents.
\end{lemme}

\begin{dem}
  On a un foncteur pleinement fidèle 
  $X\et \lra X'\et$,\ifmmode\,\fi{}  la topologie sur 
  $X\et$ est la trace de celle sur $X'\et$
  et tout objet de $X'$ peut être recouvert par des objets de $X\et$ :
  d'après le lemme de comparaison \cite[III 4.1]{SGA4}, on a donc une
  équivalence de catégories 
  $$ \cat{Sh}\,\ifmmode\,\fi{}  X\et \equiv \cat{Sh}\,\ifmmode\,\fi{}  X'\et. \qed $$
\end{dem}

\begin{theo}
  Le topos des faisceaux étales sur un champ algébrique $X$ est
  localement connexe ; on définit alors 
  $$\pi_1(X) = \pi_1 \SLC \cat{Sh}\,\ifmmode\,\fi{} X\et.$$
\end{theo}

\begin{dem}
  Voir \cite{zoonek}.
\end{dem}

\begin{rem}\label{rem:profini}%
  Ces résultats ont des analogues profinis. Soit 
  $\T$ un topos localement connexe. 

  Un objet localement constant $F$
  est dit \d{localement constant fini} si pour tout objet $U$
  trivialisant $F$,\ifmmode\,\fi{}  la restriction 
  $F|_U = (F \times U \lra U)$ est un objet constant fini de $\T/U$,\ifmmode\,\fi{} 
  \ie,\ifmmode\,\fi{}  une somme disjointe finie de copies de l'objet final de
  $\T/U$. 

  La sous-catégorie pleine $\SLCF \T$ des sommes disjointes d'objets
  localement constants finis est équivalente au topos classifiant d'un
  groupoïde profini.

  Un \d{point-base profini} de $\T$ est un morphisme de topos 
  $\cat{Ens} \lra \SLCF\T$. 

  Si $P$ est un ensemble de points-bases profinis de $\T$,\ifmmode\,\fi{}  au moins un
  dans chaque composante connexe de $\SLCF\T$,\ifmmode\,\fi{}  on définit le 
  \d{groupoïde fondamental profini} de $\T$ comme 
  $$ \wh \pi_1(\T\!,\ifmmode\,\fi{} P) = \pi_1(\SLCF \T\!,\ifmmode\,\fi{}  P). $$
  On a alors 
  $\SLCF \T \equiv \cat B \pi_1(\SLCF \T\!,\ifmmode\,\fi{}  P). $
\end{rem}

\section{Interlude sur les 2-catégories}

On peut définir le groupe fondamental (profini) d'un schéma $X$ à l'aide
de faisceaux localement constants finis ou à l'aide de revêtements
étales,\ifmmode\,\fi{}  \ie,\ifmmode\,\fi{}  de morphismes étales et finis vers~$X$. Cette seconde
description du groupe fondamental ne s'étend pas immédiatement au cas
où~$X$ est un champ algébrique : en effet,\ifmmode\,\fi{}  les morphismes étales finis
$Y \lra X$ ne forment pas une catégorie mais sont des $1$-morphismes
dans une $2$-catégorie ---
comment en faire une catégorie ?

Étant donnés une $2$-catégorie $\dcat C$ et un objet $X$ de $\dcat C$,\ifmmode\,\fi{} 
nous allons donner une définition de la $2$-catégorie 
$\dcat C / X$ des objets au dessus de $X$ et de la catégorie 
$\cat{Cat}( \dcat C/X )$ associée à cette $2$-catégorie,\ifmmode\,\fi{}  de sorte que
la notion de $2$-produit fibré au dessus de $X$ se transforme d'abord
en $2$-produit puis en produit. 
Nous utiliserons ces notions pour construire une 
$2$-catégorie $\dcat{Rev}\,\ifmmode\,\fi{} X$ puis une catégorie 
$\cat{Cat} \,\ifmmode\,\fi{}  \dcat{Rev} \,\ifmmode\,\fi{}  X$ dont les objets sont les revêtements
du champ algébrique $X$ et montrerons qu'elle est équivalente à la
catégorie 
$\LCF \cat{Sh}\,\ifmmode\,\fi{}  X$  des faisceaux étales localement constants finis sur
$X$. 

\begin{defi}
  La catégorie $\cat{Cat}\,\ifmmode\,\fi{} \dcat C$ associée à une $2$-catégorie
  $\dcat C$
  a pour objets les objets de~$\dcat C$ et pour morphismes 
  les classes d'isomorphisme de $1$-morphismes de $\dcat C$. 
\end{defi}

\begin{lemme}\label{lemme:2prodprod}%
  Un $2$-produit dans $\dcat C$ devient un produit dans
  $\cat{Cat}\,\ifmmode\,\fi{} \dcat C$. 
\end{lemme}

\begin{dem}
  Soient $X$ et $Y$ deux objets de $\dcat C$.
  Montrons que leur $2$-produit-fibré 
  $X\dtimes Y$ dans $\dcat C$ est aussi leur produit fibré dans 
  $\cat{Cat}\,\ifmmode\,\fi{} \dcat C$.

  (a)
  Le produit de $X$ et $Y$ dans 
  $\cat{Cat}\,\ifmmode\,\fi{} \dcat C$ est un diagramme 
  $$ \xymatrix{
    X \times Y \ar[r]^-{p_1} \ar[d]_{p_2} & X \\ Y,\ifmmode\,\fi{}  } $$
  tel que 
  pour tout diagramme 
  $$ \xymatrix{ T \ar[r]^{q_1} \ar[d]_{q_2} & X \\ Y } $$
  il existe un unique morphisme 
  $\phi : T \lra X \times Y$ tel que 
  $$ 
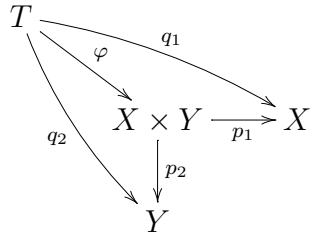
\xymatrix{ T \ar@/^/[rrd]^{q_1} \ar@/_/[rdd]_{q_2} \ar[rd]^\phi \\
    & X \times Y \ar[r]_-{p_1} \ar[d]^{p_2} & X \\ & Y } $$
  commute. 
  
  (b) D'après la définition de $\cat{Cat}\,\ifmmode\,\fi{} \dcat
  C$,\ifmmode\,\fi{}  cette propriété s'écrit aussi : pour tout
  diagramme
  $$ 
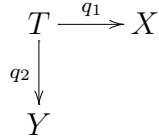
\xymatrix{ T \ar[r]^{q_1} \ar[d]_{q_2} & X \\ Y } $$
  il existe un morphisme 
  $\phi : T \lra X \times Y$ et des isomorphismes 
  \cbstart
  $\alpha_1 : p_1 \phi \lra q_1$,\ifmmode\,\fi{}  
  $\alpha_2 : p_2 \phi \lra q_2$
  $$
  \english 
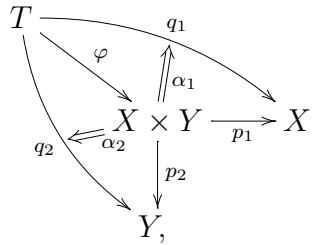
\xymatrix{ T \ar[rd]^\phi
    \ar@/^1pc/[rrd]|*{}="A" ^{q_1}
    \ar@/_1pc/[rdd]|*{}="B" _{q_2} \enlargexyentry
    A \enlargexyentry B
    \\
    & X \times Y \ar@{=>}"A"_{\alpha_1}
    \ar@{=>}"B"^{\alpha_2} \ar[r]_-{p_1}
    \ar[d]^{p_2} & X \\ & Y,\ifmmode\,\fi{}  } $$
  \cbend
  de plus,\ifmmode\,\fi{}  le
  morphisme $\phi$ est unique à isomorphisme (non
  unique,\ifmmode\,\fi{}  a priori) près.
  
  (c) Par contre la propriété $2$-universelle du
  $2$-produit-fibré s'écrit : pour tout diagramme
  $$ 
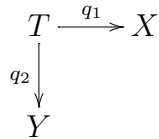
\xymatrix{ T \ar[r]^{q_1} \ar[d]_{q_2} & X \\ Y } $$
  il existe un morphisme 
  $\phi : T \lra X \times Y$ et des $2$-isomorphismes 
  \cbstart
  $\alpha_1 : p_1 \phi \lra q_1$,\ifmmode\,\fi{}  
  $\alpha_2 : p_2 \phi \lra q_2$
  \cbend
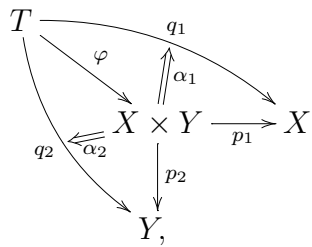
  $$ \english
  \xymatrix{ 
    T 
    \ar[rd]^\phi 
    \ar@/^1pc/[rrd]|*{}="A" ^{q_1}
    \ar@/_1pc/[rdd]|*{}="B" _{q_2}
    \enlargexyentry A
    \enlargexyentry B
    \\
    & X \times Y 
    \ar@{=>}"A"_-{\alpha_1} \ar@{=>}"B"^(.7){\alpha_2} 
    \ar[r]_-{p_1} \ar[d]^{p_2} & X \\ & Y,\ifmmode\,\fi{}  } $$
  de plus,\ifmmode\,\fi{}  le morphisme $\phi$ est unique à \emph{unique} isomorphisme
  près,\ifmmode\,\fi{}  au sens suivant : 
  si 
  $$ \english
  
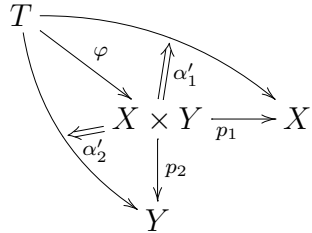
\xymatrix{ 
    T 
    \ar[rd]^\phi 
    \ar@/^1pc/[rrd]|*{}="A" 
    \ar@/_1pc/[rdd]|*{}="B" 
    \enlargexyentry A
    \enlargexyentry B
    \\
    & X \times Y 
    \ar@{=>}"A"_-{\alpha'_1} \ar@{=>}"B"^(.7){\alpha'_2}
    \ar[r]_{p_1} \ar[d]^{p_2} & X \\ & Y } $$
  est un autre diagramme comme ci-dessus,\ifmmode\,\fi{}  il existe un unique
  isomorphisme 
  $$\english  
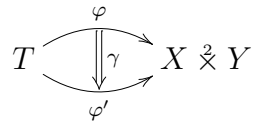
\xymatrix@!0 @C=2cm { 
    T \UN \phi{\phi'}\gamma & X \rlap{${}\dtimes Y$} } $$
  tel que le diagramme 
  $$ \english
  
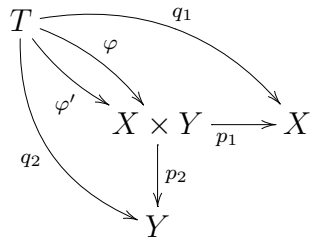
\xymatrix{ 
    T 
    \ar@/^1.5pc/[rrd]|*{}="A"  ^{q_1}
    \ar@/_1.5pc/[rdd]|*{}="B"  _{q_2}
    \ar@/^/[rd]^{\phi} \ar@/_/[rd]_{\phi'}
    \enlargexyentry A
    \enlargexyentry B
    \\
    & X \times Y 
    \ar[r]_{p_1} \ar[d]^{p_2} & X \\ & Y } $$
  soit $2$-commutatif,\ifmmode\,\fi{}  \ie,\ifmmode\,\fi{}  les diagrammes 
  \cbstart
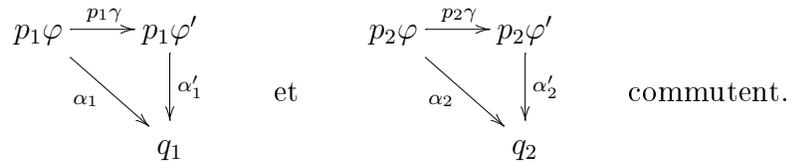
  $$ \raisebox{.5\depth}{\xymatrix{
      p_1 \phi \ar[r] ^{p_1 \gamma} \ar[rd]_{\alpha_1} &
      p_1 \phi' \ar[d]^{\alpha'_1} \\ & q_1 \strut }}
  \qquad\text{et}\qquad
  \raisebox{.5\depth}{\xymatrix{
      p_2 \phi \ar[r] ^{p_2 \gamma} \ar[rd]_{\alpha_2} &
      p_2 \phi' \ar[d]^{\alpha'_2} \\ & q_2\strut
      }}
  \qquad\text{commutent.}
  $$
  \cbend
  La condition (c) implique donc bien la condition (b). 
\end{dem}

\begin{defi}\label{defi:dcataudessus}%
  Soit $X$ un objet d'une $2$-catégorie $\dcat C$.
  La $2$-catégorie $\dcat C / X$ des objets de $\dcat C$ au dessus de
  $X$ a pour objets 
  les $1$-morphismes 
  $f : Y \lra X$ de $\dcat C$,\ifmmode\,\fi{}  pour $1$-morphismes 
  $(Y,\ifmmode\,\fi{} f) \lra (Z,\ifmmode\,\fi{}  g)$
  \cbstart
  les couples 
  $(\phi:Y\lra Z,\ifmmode\,\fi{}  \alpha: g \phi \lra f)$
  \cbend
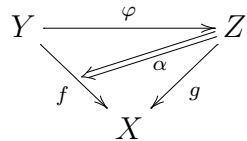
  $$ 
  \english
  \xymatrix{
    Y \ar[rr]^\phi \ar[rd]|*{}="A"_f 
    \enlargexyentry A
    && Z \ar[ld]^{g} \ar@{=>}"A"^\alpha \\ & X } $$
  et pour $2$-morphismes 
  $(\phi,\ifmmode\,\fi{}  \alpha) \lra (\psi,\ifmmode\,\fi{}  \beta)$
  les $2$-morphismes
  $$ \english 
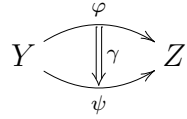
\xymatrix @!0 @C=2cm { Y \UN \phi\psi\gamma & Z } $$
  tels que le diagramme 
  \cbstart
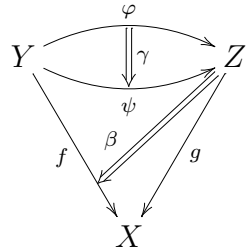
  $$\english
  \xymatrix{
    Y 
    \UN[rr]{\phi}{\psi}{\gamma}
    \ar[rdd]|(.7)*{}="Z"_f && Z \ar[ldd]^g \\\\ & X 
    \ar@{=>}[ruu];"Z"_(.8)\beta 
    } $$
  soit $2$-commutatif,\ifmmode\,\fi{}  \ie,\ifmmode\,\fi{}  le diagramme suivant
  commute. 
  $$ 
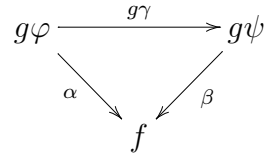
\xymatrix{
    g \phi \ar[rr]^{ g \gamma } \ar[rd]_\alpha && 
    g \psi \ar[ld]^\beta \\ & f } $$
  \cbend
\end{defi}

\begin{lemme}\label{lemme:2prod2prodfib}%
  Le foncteur d'oubli 
  $$\definefunction{\dcat C/X}{\dcat C}{(Y\lra X)}{Y}$$
  transforme les $2$-produits en $2$-produits-fibrés. 
\end{lemme}

\begin{dem}
  (a) Un $2$-produit 
  $Y \times Z$ dans $\dcat C / X$ est un diagramme 
  $$ \xymatrix{ Y \times Z \ar[r]^-{p_2} \ar[d]_{p_1} & Z \\ Y } $$
  tel que pour tout diagramme 
  $$ \xymatrix{ T \ar[r]^{q_2} \ar[d]_{q_1} & Z \\ Y } $$
  il existe un morphisme $\phi$ et des $2$-morphismes $\alpha_1$,\ifmmode\,\fi{}  $\alpha_2$,\ifmmode\,\fi{} 
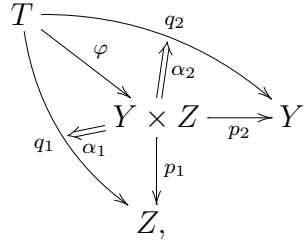
  $$ \english
  \xymatrix{ 
    T 
    \ar@/^1pc/[rrd]|*{}="A"^{q_2} 
    \ar@/_1pc/[rdd]|*{}="B"_{q_1}
    \ar[rd]^\phi 
    \enlargexyentry A
    \enlargexyentry B
    \\ & Y \times Z
    \ar@{=>}"A"_-{\alpha_2} \ar@{=>}"B"^(.7){\alpha_1\vphantom{b}}
    \ar[r]_-{p_2} \ar[d]^{p_1} & Y \\ & Z,\ifmmode\,\fi{}  } $$
  tels que pour tout autre diagramme 
  $$ \english
  
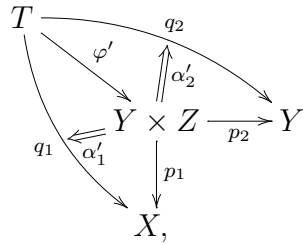
\xymatrix{ 
    T 
    \ar@/^1pc/[rrd]|*{}="A" ^{q_2}
    \ar@/_1pc/[rdd]|*{}="B" _{q_1}
    \ar[rd]^{\phi'} 
    \enlargexyentry A
    \enlargexyentry B
    \\
    & Y \times Z
    \ar@{=>}"A"_-{\alpha'_2} \ar@{=>}"B"^(.7){\alpha'_1}
    \ar[r]_-{p_2} \ar[d]^{p_1} & Y \\ & X,\ifmmode\,\fi{}  } $$
  il existe un unique isomorphisme 
  $$ \english 
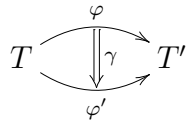
\xymatrix @!0 @C=2cm { T \UN{\phi}{\phi'}\gamma & T' }$$
  tel que le diagramme 
  $$ \english
  
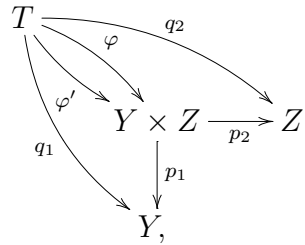
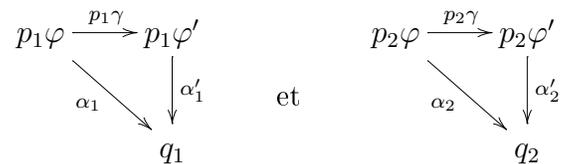
\xymatrix{ 
    T 
    \ar@/^1pc/[rrd]|*{}="A" ^{q_2}
    \ar@/_1pc/[rdd]|*{}="B" _{q_1}
    \ar@/^/[rd]^\phi 
    \ar@/_/[rd]_{\phi'}
    \enlargexyentry A
    \enlargexyentry B 
    \\
    & Y \times Z
    \ar[r]_-{p_2} \ar[d]^{p_1} & Z \\ & Y,\ifmmode\,\fi{}  } $$
  soit $2$-commutatif,\ifmmode\,\fi{}  \ie,\ifmmode\,\fi{}  
  $$ \raisebox{.5\depth}{\xymatrix{
    p_1 \phi \ar[r] ^{ p_1 \gamma } \ar[rd]_{ \alpha_1 } &
    p_1 \phi' \ar[d]^{\alpha'_1} \\ & q_1 \strut }}
  \qquad\text{et}\qquad
  \raisebox{.5\depth}{\xymatrix{
    p_2 \phi \ar[r] ^{ p_2 \gamma } \ar[rd]_{ \alpha_2 } &
    p_2 \phi' \ar[d]^{\alpha'_2} \\ & q_2 \strut }}
  $$
  
  (b) D'après la définition de 
  $\dcat C / X$,\ifmmode\,\fi{}  le $2$-produit $X \dtimes Y$ est un diagramme 
  $$ \english 
  
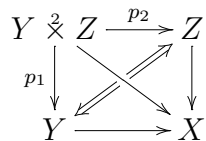
\xymatrix{
    Y \dtimes Z \ar[r]^-{p_2}  \ar[d]_{p_1} \ar[rd]|*{}="A" 
    \enlargexyentry A
    & Z \ar[d] \ar@{=>}"A";[] \\
    Y \ar[r] \ar@{=>}"A";[] & X } $$
  tel que pour tout diagramme 
  $$ \english 
  \xymatrix{
    T \ar[r]^{q_2} \ar[d]_{q_1} \ar[rd]|*{}="A" 
    \enlargexyentry A
    & Z \ar[d] \ar@{=>}"A";[] \\
    Y \ar[r] \ar@{=>}"A";[] & X } $$
  il existe $\phi$,\ifmmode\,\fi{}  $\alpha_1$,\ifmmode\,\fi{}  $\alpha_2$ rendant le diagramme
  suivant $2$-commutatif,\ifmmode\,\fi{} 
  $$ \english
  \xymatrix{ 
    T 
    \ar@/^1pc/[rrd]|*{}="A" ^{q_2}
    \ar@/_1pc/[rdd]|*{}="B" _{q_1}
    \ar[rd]^{\phi} 
    \enlargexyentry A
    \enlargexyentry B
    \\
    & Y \dtimes Z
    \ar[rd]|*{}="C"
    \ar@{=>}"A"_-{\alpha_2} \ar@{=>}"B"^(.7){\alpha_1\vphantom{b}}
    \ar[r]_{p_2} \ar[d]^{p_1} 
    \enlargexyentry C
    & Z \ar[d] \ar@{=>}"C";[]
    \\ & Y \ar[r] \ar@{=>}"C";[] & X } $$
  et tel que pour tout autre diagramme $2$-commutatif
  $$ \english
  \xymatrix{ 
    T 
    \ar@/^1pc/[rrd]|*{}="A" ^{q_2}
    \ar@/_1pc/[rdd]|*{}="B" _{q_1}
    \ar[rd]^{\phi'} 
    \enlargexyentry A
    \enlargexyentry B
    \\
    & Y \dtimes Z
    \ar[rd]|*{}="C"
    \ar@{=>}"A"_-{\alpha'_2} \ar@{=>}"B"^(.7){\alpha'_1}
    \ar[r]_{p_2} \ar[d]^{p_1} 
    \enlargexyentry C
    & Z \ar[d] \ar@{=>}"C";[]
    \\ & Y \ar[r] \ar@{=>}"C";[] & X } $$
  il existe un unique isomorphisme 
  $$ \english \xymatrix @!0 @C=2cm{ 
    T \UN{\phi}{\phi'}\gamma & Y \rlap{${}\dtimes Z$} } $$
  tel que 
  $$ \raisebox{.5\depth}{\xymatrix{
    T \ar@/^/[rr]^\phi \ar@/_/[rr]_{\phi'} \ar[rd] && 
    Y \dtimes Z \ar[ld] \\ & X }}
  \qquad\text{et}\qquad
  \english
  \raisebox{.5\depth}{\xymatrix{
    T 
    \ar@/^1pc/[rrd]|*{}="A" ^{q_2}
    \ar@/_1pc/[rdd]|*{}="B" _{q_1}
    \ar@/^/[rd]^{\phi}
    \ar@/_/[rd]_{\phi'}
    \enlargexyentry A
    \enlargexyentry B
    \\
    & Y \dtimes Z
    \ar[rd]|*{}="C"
    \ar[r]_{p_2} \ar[d]^{p_1}     
    \enlargexyentry C
    & Z \ar[d] 
    \\ & Y \ar[r] 
    & X }}
  $$
  soient $2$-commutatifs (mais comme le premier diagramme est contenu
  dans le second,\ifmmode\,\fi{}  on peut l'oublier).
    
  (c) En remarquant que tout diagramme 
  $$ \english
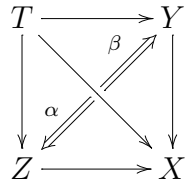
\xymatrix@!0 @R=2cm @C=2cm{ 
    T \ar[r] \ar[d] \ar[rd]|*{}="A" 
    \enlargexyentry A
    & Y \ar[d] \ar@{=>}"A";[]^-\beta\\
    Z \ar[r] \ar@{=>}"A";[]_-\alpha & X } $$
  définit un diagramme 
  $$ \english
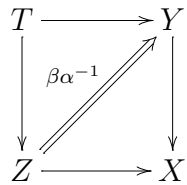
\xymatrix@!0 @R=2cm @C=2cm{  
    T \ar[r] \ar[d] & Y \ar[d] \\
    Z \ar[r] \ar@{=>}[ru]^{\beta\alpha^{-1}\!\!\!} & X } $$
  et que réciproquement,\ifmmode\,\fi{}  tout diagramme 
  $$ \english
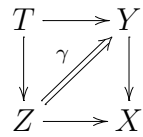
\xymatrix{
    T \ar[r] \ar[d] & Y \ar[d] \\
    Z \ar[r] \ar@{=>}[ru]^\gamma & X } $$  
  définit un diagramme 
  $$ \english
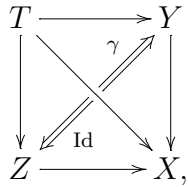
\xymatrix@!0 @R=2cm @C=2cm{ 
    T \ar[r] \ar[d] \ar[rd]|*{}="A" 
    \enlargexyentry A
    & Y \ar[d] \ar@{=>}"A";[]^-{\gamma}\\
    Z \ar[r] \ar@{=>}"A";[]^-{\id} & X, } $$
  on voit que le produit fibré 
  $Y \dtimes Z$ dans $\dcat C/X$ est un diagramme 
  $$ \english
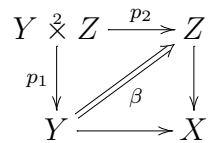
\xymatrix{ 
    Y \dtimes Z \ar[r]^-{p_2} \ar[d]_{p_1} & Z \ar[d] \\
    Y \ar[r] \ar@{=>}[ru]_\beta & X } $$
  tel que pour tout 
  $$ \english
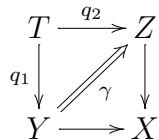
\xymatrix{ 
    T \ar[r]^-{q_2} \ar[d]_{q_1} & Z \ar[d] \\
    Y \ar[r] \ar@{=>}[ru]_\gamma & X } $$
  il existe un diagramme $2$-commutatif 
  $$ \english
  
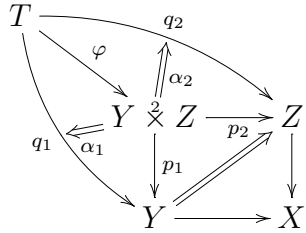
\xymatrix{ 
    T 
    \ar@/^1pc/[rrd]|*{}="A" ^{q_2}
    \ar@/_1pc/[rdd]|*{}="B" _{q_1}
    \ar[rd]^{\phi}
    \enlargexyentry A
    \enlargexyentry B
    \\
    & Y \dtimes Z
    \ar@{=>}"A"_{\alpha_2} \ar@{=>}"B"^(.7){\alpha_1\vphantom{b}}
    \ar[r]_-{p_2} \ar[d]^{p_1} & Z \ar[d] 
    \\ & Y \ar[r] \ar@{=>}[ru] & X } $$
  et tel que pour tout autre diagramme $2$-commutatif
  $$ \english
  
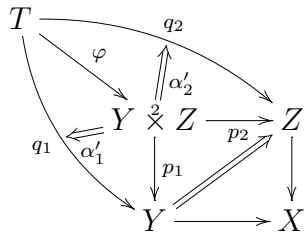
\xymatrix{ 
    T 
    \ar@/^1pc/[rrd]|*{}="A" ^{q_2}
    \ar@/_1pc/[rdd]|*{}="B" _{q_1}
    \ar[rd]^{\phi} 
    \enlargexyentry A
    \enlargexyentry B
    \\
    & Y \dtimes Z
    \ar@{=>}"A"_{\alpha'_2} \ar@{=>}"B"^(.7){\alpha'_1\vphantom{b}}
    \ar[r]_-{p_2} \ar[d]^{p_1} & Z \ar[d] 
    \\ & Y \ar[r] \ar@{=>}[ru] & X } $$
  il existe un unique isomorphisme 
  $$ \english
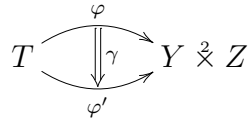
\xymatrix@!0 @C=2cm{ 
    T \UN{\phi}{\phi'}\gamma & Y \rlap{\ensuremath{{}\dtimes Z}} } $$
  tel que 
  $$ \english
  
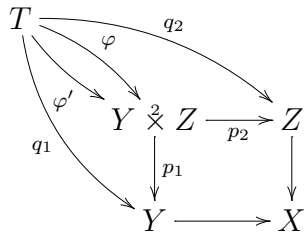
\xymatrix{ 
    T 
    \ar@/^1pc/[rrd]|*{}="A" ^{q_2}
    \ar@/_1pc/[rdd]|*{}="B" _{q_1}
    \ar@/^/[rd]^\phi
    \ar@/_/[rd]_{\phi'}
    \enlargexyentry A
    \enlargexyentry B
    \\
    & Y \dtimes Z
    \ar[r]_-{p_2} \ar[d]^{p_1}     
    & Z \ar[d] 
    \\ & Y \ar[r] 
    & X } $$
  soit $2$-commutatif.
  Mais il s'agit de la propriété $2$-universelle définissant le
  $2$-produit-fibré 
  $Y \dtimes_X Z$ dans $\dcat C$. 
\end{dem}

\begin{rem}
  En général,\ifmmode\,\fi{}  
  $\cat{Cat}\,\ifmmode\,\fi{} (\dcat C/X) \not\equiv (\cat{Cat}\,\ifmmode\,\fi{} \dcat C)/X$.
  Considérons par exemple la $2$-catégorie $\dcat C$ définie par 
  \begin{align*}
    \ob\dcat C &= \{\,\ifmmode\,\fi{} X,\ifmmode\,\fi{}  Y\,\ifmmode\,\fi{} \} 
    &
    \ob\cat{hom}(Y,\ifmmode\,\fi{}  X) &= \{a\}    
    \\
    \cat{hom}(X,\ifmmode\,\fi{}  X) &= \{\id_X\} 
    &
    \hom(a,\ifmmode\,\fi{} a) &= \{\id_a,\ifmmode\,\fi{}  \alpha\}
    \\
    \cat{hom}(Y,\ifmmode\,\fi{} Y) &= \{\id_Y\} 
    &
    \alpha^2 &= \id_a 
    \\
    \cat{Hom}( X,\ifmmode\,\fi{} Y ) &= \nothing.
  \end{align*}
  $$ \includegraphics{1.1} $$
  La catégorie $\cat{Cat}\,\ifmmode\,\fi{} (\dcat C /X)$ a alors 
  $a$ pour seul objet et $\{\id_a,\ifmmode\,\fi{}  \alpha\}$ pour morphismes,\ifmmode\,\fi{}  alors
  que $(\cat{Cat}\,\ifmmode\,\fi{} \dcat C)/X$ 
  a $a$ pour seul objet et $\id_a$ pour seul
  morphisme. 
\end{rem}

Nous n'utiliserons pas le résultat suivant par la suite. 

\begin{lemme}
  Soit $U \lra X$ un objet d'une $2$-catégorie $\dcat C$.
  On a une équivalence de $2$-catégories 
  $(\dcat C/X)/U \equiv \dcat C / U$.  
\end{lemme}
\begin{dem}
  Notons $u : U \lra X$. Nous allons décrire les $2$-catégories 
  $\dcat C / U$ et $\dcat C / X / U$ et constater qu'elles sont
  équivalentes. 

  (a)
  La catégorie $\dcat C / U$ a pour objets les 
  $f' : F \lra U$,\ifmmode\,\fi{}  pour morphismes 
  $(F,\ifmmode\,\fi{}  f') \lra (G,\ifmmode\,\fi{} g')$ les 
  \cbstart
  $(h:F\lra G,\ifmmode\,\fi{} \delta:g'h\lra f')$
  \cbend
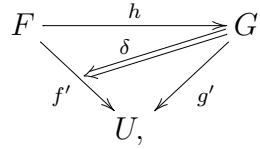
  $$ 
  \english
  \xymatrix{
    F 
    \ar[rr]^{h} 
    \ar[rd]|*{}="A"_{f'}
    \enlargexyentry A
    && 
    G
    \ar[ld]^{g'} 
    \ar@{=>}"A"_(.7)\delta
    \\ & 
    U,\ifmmode\,\fi{} 
    }
  $$
  et les $2$-morphismes 
  $$ 
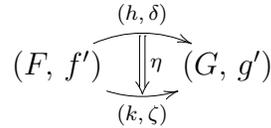
\xymatrix{
    (F,\ifmmode\,\fi{} f') \UN{(h,\ifmmode\,\fi{} \delta)}{(k,\ifmmode\,\fi{} \zeta)}{\eta}
    & (G,\ifmmode\,\fi{} g') }
  $$
  sont les 2-morphismes 
  $$
    
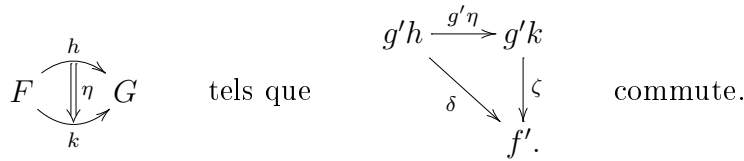
\xymatrix{
      F \UN hk\eta & G 
      }
  \qquad\text{tels que}\qquad
  \raisebox{.5\depth}{\xymatrix{
      g' h \ar[r]^{g' \eta} \ar[rd]_\delta & g' k
      \ar[d]^\zeta \\ & f'. }}
  \qquad\text{commute.}
  $$

  (b)
  La catégorie $\dcat C / X / U$ a pour objets les 
  \cbstart
  $(F,\ifmmode\,\fi{}  f:F\lra X,\ifmmode\,\fi{}  f':F\lra U,\ifmmode\,\fi{}  \alpha:uf'\lra f)$
  \cbend
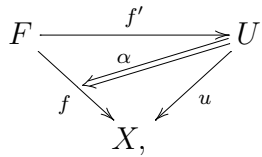
  $$ 
  \english
  \xymatrix{
    F 
    \ar[rr]^{f'} 
    \ar[rd]|*{}="A"_f
    \enlargexyentry A
    && 
    U 
    \ar[ld]^u 
    \ar@{=>}"A"_(.7)\alpha
    \\ & 
    X,\ifmmode\,\fi{} 
    }
  $$
  pour $1$-morphismes 
  $(F,\ifmmode\,\fi{}  f,\ifmmode\,\fi{} f',\ifmmode\,\fi{} \alpha) \lra (G,\ifmmode\,\fi{}  g,\ifmmode\,\fi{}  g',\ifmmode\,\fi{}  \beta)$ 
  \cbstart
  les $(h:F\lra G,\ifmmode\,\fi{}  \gamma:gh\lra f,\ifmmode\,\fi{}  \delta:g'h\lra
  f')$,\ifmmode\,\fi{} 
  \cbend
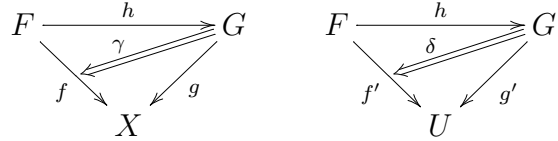
  $$
  \english
  \xymatrix{
    F 
    \ar[rr]^{h} 
    \ar[rd]|*{}="A"_f
    \enlargexyentry A
    && 
    G 
    \ar[ld]^g
    \ar@{=>}"A"_(.7)\gamma
    \\ & 
    X
    }
  \qquad
  \xymatrix{
    F 
    \ar[rr]^{h} 
    \ar[rd]|*{}="A"_{f'}
    \enlargexyentry A
    && 
    G 
    \ar[ld]^{g'}
    \ar@{=>}"A"_(.7)\delta
    \\ & 
    U
    }
  $$
  tels que  le diagramme 
  $$
  
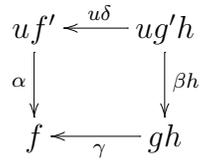
\xymatrix{
    u f' \ar[d]_\alpha & u g' h \ar[l]_{u \delta} \ar[d]^{\beta h} \\ 
    f & g h \ar[l]^\gamma 
    }
  $$
  commute,\ifmmode\,\fi{}  
  (on remarquera que la condition définit $\gamma$ en fonction de
  $\delta$,\ifmmode\,\fi{}  on peut donc voir un $1$-morphisme comme la donnée de 
  $h$ et $\delta$ tout seuls)
  et pour $2$-morphismes
  $$ 
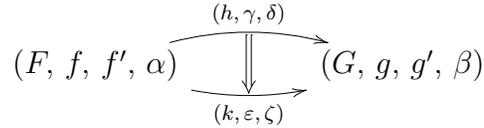
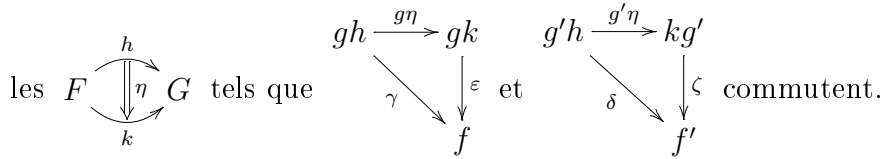
\xymatrix{
    {(F,\ifmmode\,\fi{} f,\ifmmode\,\fi{} f',\ifmmode\,\fi{} \alpha)\quad}
    \UN{(h,\ifmmode\,\fi{} \gamma,\ifmmode\,\fi{} \delta)}{(k,\ifmmode\,\fi{} \epsilon,\ifmmode\,\fi{} \zeta)}{}
    & 
    {\quad(G,\ifmmode\,\fi{} g,\ifmmode\,\fi{} g',\ifmmode\,\fi{} \beta)}
    }
  $$
  les 
  $ \xymatrix{
    F \UN hk\eta & G 
    }
  $ 
  tels que 
  $
  \raisebox{.5\depth}{\xymatrix{ 
      g h \ar[r]^{g \eta} \ar[rd]_\gamma & g k
      \ar[d]^\epsilon \\ & f }}
  $
  et
  $
  \raisebox{.5\depth}{\xymatrix{
      g' h \ar[r]^{g' \eta} \ar[rd]_\delta & k g' \ar[d]^\zeta \\ &
      f' }}
  $
  commutent.
  
  (c) On définit alors un 2-foncteur 
  $$ \english 
  \left\{
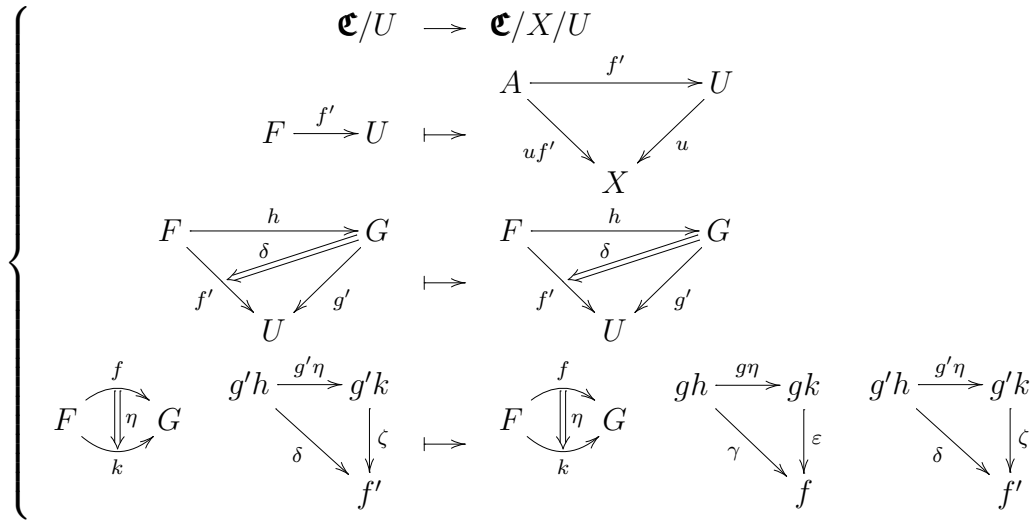
    \begin{array}{rcl}
      \dcat C / U & \lra & \dcat C / X / U \\
      \raisebox{.5\depth}{\xymatrix{F \ar[r]^{f'} & U}}
      & \lmt &
      \raisebox{.5\depth}{\xymatrix{
          A \ar[rr]^{f'} \ar[rd]_{u f'} && U \ar[ld]^u \\ & X}}
      \\
      \raisebox{.5\depth}{\xymatrix{
        F 
        \ar[rr]^{h} 
        \ar[rd]|*{}="A"_{f'}
        \enlargexyentry A
        && 
        G 
        \ar[ld]^{g'}
        \ar@{=>}"A"_(0.7)\delta
        \\ & 
        U
        }}
      & \lmt &
      \raisebox{.5\depth}{\xymatrix{
        F 
        \ar[rr]^{h} 
        \ar[rd]|*{}="A"_{f'}
        \enlargexyentry A
        && 
        G 
        \ar[ld]^{g'}
        \ar@{=>}"A"_(0.7)\delta
        \\ & 
        U
        }}
      \\
      \raisebox{.5\depth}{\xymatrix{ F \UN fk\eta & G }}
      \quad
      \raisebox{.5\depth}{\xymatrix{ 
          g' h \ar[r]^{g' \eta} \ar[rd]_\delta & g' k \ar[d]^\zeta \\ &
          f' }}
      & \lmt &
      \raisebox{.5\depth}{\xymatrix{ F \UN fk\eta & G }}
      \quad
      \raisebox{.5\depth}{\xymatrix{ gh \ar[r]^{g\eta} \ar[rd]_\gamma
          & gk \ar[d]^\epsilon
          \\ & f }}
      \quad
      \raisebox{.5\depth}{\xymatrix{ 
          g' h \ar[r]^{g' \eta} \ar[rd]_\delta & g' k \ar[d]^\zeta \\ &
          f' }}
    \end{array}
  \right.
  $$

  Ce 2-foncteur est essentiellement surjectif,\ifmmode\,\fi{}  car pour tout objet 
  $(F,\ifmmode\,\fi{} f,\ifmmode\,\fi{} f',\ifmmode\,\fi{} \alpha)$ de $\dcat C / X / U$,\ifmmode\,\fi{}  on a un isomorphisme 
  $$ \english \xymatrix {
    (F,\ifmmode\,\fi{}  u f',\ifmmode\,\fi{}  f',\ifmmode\,\fi{}  \id) 
    \ar[rr]^{(\id,\ifmmode\,\fi{}  \id,\ifmmode\,\fi{}  \alpha)} &&
    (F,\ifmmode\,\fi{} f,\ifmmode\,\fi{} f',\ifmmode\,\fi{} \alpha)
    }
  $$

  Montrons que ce 2-foncteur est pleinement fidèle,\ifmmode\,\fi{}  \ie,\ifmmode\,\fi{}  que l'on a des
  équivalences (et même,\ifmmode\,\fi{}  en fait,\ifmmode\,\fi{}  des isomorphismes) de catégories 
  $$ \cat{hom}_{\dcat C/U}\bigl(
  (F,\ifmmode\,\fi{} f'),\ifmmode\,\fi{}  (G,\ifmmode\,\fi{} g') 
  \bigr) \lra \cat{hom}\bigl(
  (F,\ifmmode\,\fi{}  uf',\ifmmode\,\fi{}  f',\ifmmode\,\fi{}  \id),\ifmmode\,\fi{}  
  (G,\ifmmode\,\fi{}  ug',\ifmmode\,\fi{}  g',\ifmmode\,\fi{}  \id)
  \bigr). $$
  Tout d'abord les objets,\ifmmode\,\fi{}  à droite comme à
  gauche,\ifmmode\,\fi{}  sont les 
  diagrammes
  $$ \xymatrix{
    F 
    \ar[rr]^{h} 
    \ar[rd]|*{}="A"_{f'}
    \enlargexyentry A
    && 
    G 
    \ar[ld]^{g'}
    \ar@{=>}"A"_(.7)\delta
    \\ & 
    U
    }
  $$
  D'autre part,\ifmmode\,\fi{}  à gauche les morphismes sont les 
  $$ 
  \raisebox{.5\depth}{\xymatrix{ F \UN hk\eta & G }}
  \qquad 
  \raisebox{.5\depth}{\xymatrix{ 
      g' h \ar[r]^{g' \eta} \ar[rd]_\delta & g' k \ar[d]^\zeta
      \\ & f' }}
  $$
  et à droite ce sont les 
  $$ 
  \raisebox{.5\depth}{\xymatrix{ F \UN hk\eta & G }}
  \qquad 
  \raisebox{.5\depth}{\xymatrix{ 
      g' h \ar[r]^{g' \eta} \ar[rd]_\delta & g' k \ar[d]^\zeta
      \\ & f' }}
  \qquad 
  \text{tels que}
  \qquad
  \raisebox{.5\depth}{\xymatrix{ 
      gh \ar[r]^{g\eta}\ar[rd]_\gamma & gk \ar[d]^\epsilon \\ &
      f }}
  \text{commute}
  $$
  or cette dernière condition est automatiquement vérifiée. 
\end{dem}

\section{Interprétation géométrique des revêtements d'un champ
  algébrique}

On note $\dcat{Champs}$ la $2$-catégorie des champs algébriques. 
Si $X$ est un champ algébrique,\ifmmode\,\fi{}  on note 
$\LCF\cat{Sh}\,\ifmmode\,\fi{} X$ la catégorie des faisceaux étales localement
constants finis sur $X$ : on sait que cette catégorie est galoisienne
(theorème~\ref{theo:pi1topos} et remarque~\ref{rem:profini}) 
et qu'elle permet de calculer le groupe fondamental profini de
$X$. Nous allons montrer qu'elle est équivalente à une catégorie dont
les objets sont les revêtements étales de $X$. 

\begin{defi}\label{def:revchamp}
  La $2$-catégorie $\dcat{Rev}\,\ifmmode\,\fi{} X$ des revêtements 
  d'un champ algébrique $X$ est la sous-$2$-catégorie pleine 
  de $\dcat{Champs}/X$ constituée des morphismes 
  $Y \lra X$ étales finis. 
\end{defi}

\begin{theo}\label{theo:revchamp}%
  Soit $X$ un champ algébrique.
  On a une équivalence de catégories 
  $$ \LCF \cat{Sh}\,\ifmmode\,\fi{}  X \equiv \cat{Cat}\,\ifmmode\,\fi{}  \dcat{Rev}\,\ifmmode\,\fi{}  X. $$
\end{theo}

\begin{dem}
  D'après \cite{vistoli}, on peut supposer que $X = [U/R]$, 
  où $R \dlra U$ est un groupoïde étale dont la diagonale est
  quasi-compacte et séparée et on peut identifier les faisceaux
  étales sur $X$ aux faisceaux équivariants sur ce groupoïde. Nous
  regarderons donc les objets de $\LCF \cat{Sh}\,\ifmmode\,\fi{}  X$ comme des
  revêtements étales équivariants du groupoïde 
  $R \dlra U$. 

  Considérons alors le foncteur 
  $$ \f F : \definefunction{\LCF\cat{Sh}\,\ifmmode\,\fi{} X}{\cat{Cat}\,\ifmmode\,\fi{} \dcat{Rev}\,\ifmmode\,\fi{} X}%
  {\raisebox{.5\depth}{\xymatrix{&F\ar[d]\\R\dar[r]&U}}}{
    \raisebox{.5\depth}{\xymatrix{[F/R]\ar[d]\\ X.}}}$$

  (a) Montrons que le foncteur $\f F$ est bien défini.

  (a1) Le groupoïde $R \times F \dlra F$ est étale.

  (a2) Montrons que sa diagonale $R \times F \dlra F$ est séparée : on
  a un carré commutatif (mais pas cartésien) 
  $$ \english \xymatrix @!0 @R=1.5cm @C=3cm {
    R \Times_U F \ar[r]^{ (\pr_2,\ifmmode\,\fi{}  \alpha) } \ar[d]_{1\times f} &
    F \times F \ar[d]^{f \times f} \\ 
    R \ar[r]_{(s,\ifmmode\,\fi{} b)} & U \times U,\ifmmode\,\fi{} } $$
  où $1 \times f$ et $(s,\ifmmode\,\fi{} b)$ sont séparés,\ifmmode\,\fi{}  donc 
  $(f \times f) \circ (\pr_2,\ifmmode\,\fi{}  \alpha) = 
  (s,\ifmmode\,\fi{} b) \circ (1 \times f)$ est séparé,\ifmmode\,\fi{}  donc 
  $(\pr_2,\ifmmode\,\fi{}  \alpha)$ est séparé. 

  (a3) Montrons que la diagonale $R \times F \dlra F$ est
  quasi-compacte. Rappelons qu'un morphisme de schémas 
  $\phi : A \lra B$ est quasi-compact \ssi il existe un recouvrement
  de $B$ par des ouverts affines $U_i$ et,\ifmmode\,\fi{}  pour tout $i$,\ifmmode\,\fi{}  
  un nombre fini d'ouverts
  affines $(V_{i,\ifmmode\,\fi{} j})_{j\in J_i}$ 
  de $A$ qui recouvrent $\phi^{-1}(U_i)$,\ifmmode\,\fi{}  \ie,\ifmmode\,\fi{}  
  $\bigunion_{j \in J_i} V_{i,\ifmmode\,\fi{} j} \supset \phi^{-1}(U_i)$ ; c'est alors
  vrai pour n'importe quel recouvrement $(U_i)_{i\in I}$ par des
  ouverts affines. 
  On a un carré commutatif 
  $$ \english \xymatrix @!0 @R=1.5cm @C=3cm {
    R \Times_U F \ar[r]^{ (\pr_2,\ifmmode\,\fi{}  \alpha) } \ar[d]_{1\times f} &
    F \times F \ar[d]^{f \times f} \\ 
    R \ar[r]_{(s,\ifmmode\,\fi{} b)} & U \times U} $$
  dont les morphismes verticaux sont des revêtements étales et le
  morphisme du bas est quasi-compact. Soit $(V_i)_{i\in I}$ un
  recouvrement de $U \times U$ par des ouverts affines et,\ifmmode\,\fi{}  pour chaque
  $i\in I$,\ifmmode\,\fi{}  
  un recouvrement $(R_{i,\ifmmode\,\fi{} j})_{j \in J_i}$ de
  $(s,\ifmmode\,\fi{} b)^{-1}(V_i)$ par un nombre fini d'ouverts affines. 
  Comme $1 \times f$ est quasi-compact,\ifmmode\,\fi{}  on peut recouvrir chaque 
  $(1 \times f)^{-1} (R_{ij})$ par un nombre fini d'ouverts affines 
  $(W_{ijk})_{k \in K_{ij}}$. D'autre part,\ifmmode\,\fi{}  comme 
  $f \times f : F \times F \lra U \times U$ est quasi-compact,\ifmmode\,\fi{}  on peut
  recouvrir  chaque 
  $(f\times f)^{-1}(V_i)$ 
  par un nombre fini d'ouverts affines 
  $(G_{il})_{l\in L_i}$. Dès lors,\ifmmode\,\fi{}  les 
  $(W_{ijk})_{j \in J_i,\ifmmode\,\fi{} \ k \in K_{ij}}$ sont des ouverts affines qui
  recouvrent chaque 
  $(\pr_2,\ifmmode\,\fi{}  \alpha)^{-1}( G_{il} )$.

  (a4) D'après le lemme \ref{lemme:prescarcart},\ifmmode\,\fi{}  on a un carré
  $2$-cartésien 
  $$ \english \xymatrix{ 
    F \ar@{->>}[r]^-{\text{étale}} 
    \ar[d]_{\text{étale}} 
    \cartesien &
    \relax[F/R] \ar[d] \\
    U \ar@{->>}[r]_-{\text{étale}} &
    \relax[U/R]\rlap{,\ifmmode\,\fi{} }} $$
  donc le morphisme 
  $ [F/R] \lra [U/R]$ est étale,\ifmmode\,\fi{}  quasi-fini et propre 
  \cite[4.11]{DM}.

  (b) Montrons que le foncteur $\f F$ est essentiellement surjectif.
  Soit $Y \lra X$ un morphisme étale et fini. Le morphisme 
  $Y \times_X U \lra U$ est toujours étale et fini et on a un
  morphisme de groupoïdes 
  $$ \xymatrix{
    (Y \Times_X U) \Times_Y (Y \Times_X U) \ar[d]
    \dar[r] & (Y \Times_X U)  \cartesien \ar[r] \ar[d] 
    & Y \ar[d] \\
    U \Times_X U \dar[r] & U \ar[r] & X. } $$
  Mais comme 
  \begin{align*}
    (Y \Times_X U) \Times_Y (Y \Times_X U)
    & \iso U \Times_X (Y \Times_X U) \\
    & \iso (U \Times_X U) \Times_U (Y \Times_X U),\ifmmode\,\fi{} 
  \end{align*}
  on voit que la présentation de $Y$ ainsi obtenue est de la forme
  requise : c'est un revêtement équivariant
  de $R \dlra U$.

  (c) Pour montrer que le foncteur $\f F$ est pleinement fidèle,\ifmmode\,\fi{} 
  soient 
  $F$ et $G$ des revêtements équivariants du groupoïde 
  $R \dlra U$.
  Montrons que l'application 
  $$ \hom_{ R \dlra U }( F,\ifmmode\,\fi{}  G ) \lra \cat{hom}_{[U/R]} 
  \bigl( [F/R],\ifmmode\,\fi{}  [G/R] \bigr) / \text{isomorphismes}$$ 
  est injective. 
  Soient 
  $\phi,\ifmmode\,\fi{} \,\ifmmode\,\fi{}  \psi : F \lra G$ deux morphismes équivariants et 
  $$ \xymatrix{
    \relax[F/R] \UUN{\phi/R}{\psi/R}{\alpha} & \relax[G/R] } $$
  un isomorphisme entre leurs images.
  Si l'on prend la fibre de ces catégories fibrées au dessus de~$F$,\ifmmode\,\fi{} 
  on obtient un isomorphisme 
  $$ \english \xymatrix @!0 @C=8pc {
    \relax[F/R]_F \UUN{(\phi/R)_F}{(\psi/R)_F}{\alpha_F} & \relax[G/R]_F. } $$
  D'après la description des fibres d'un champ algébrique quotient
  \cite{vistoli}, cela nous donne un isomorphisme 
  $$
  \english
  \savebox{\tempboxa}{\ensuremath{\dlimind}}
  \savebox{\tempboxb}{%
    \raisebox{.5\depth}{\xymatrix{
        T \Times_F T \dar[d] \\ T }}}
  \savebox{\tempboxc}{%
    \raisebox{.5\depth}{\xymatrix{
        R \Times_U F  \dar[d] \\ F }}}
  \savebox{\tempboxd}{%
    \raisebox{.5\depth}{\xymatrix{
        T \Times_F T \dar[d] \\ T }}}
  \savebox{\tempboxe}{%
    \raisebox{.5\depth}{\xymatrix{
        R \Times_U G  \dar[d] \\ G }}}
  \xymatrix @!0 @C=3cm {
    \hbox to 1cm{\hss\ensuremath{
        \mathop{\usebox{\tempboxa}}\limits_{
          \substack{
            T \lra F \\
            \text{surjectif et étale}
            }}
        \cat{hom}\left(\usebox{\tempboxb},\ifmmode\,\fi{} \usebox{\tempboxc}\right)
        }}
    \UN{(\phi/R)_F}{(\psi/R)_F}{\alpha_F} & 
    \hbox to 1cm{\ensuremath{
        \mathop{\usebox{\tempboxa}}\limits_{
          \substack{
            T \lra F \\
            \text{surjectif et étale}
            }}
        \cat{hom}\left(\usebox{\tempboxd},\ifmmode\,\fi{} \usebox{\tempboxe}\right)
        }\hss}
    } $$
  La catégorie de gauche possède un objet particulier,\ifmmode\,\fi{}  correspondant
  au morphisme canonique $F \lra [F/R]$ : regardons son image $\gamma$
  par
  $\alpha_F$. C'est un 2-morphisme de groupoïdes 
  $$ \english\xymatrix{
    T \Times_F T \ar[d] \dar[r] & T \ar[d] 
    \ar[ddl]
    |!{[ld];[d]}\hole
    ^(.3)\gamma \\ 
    R \Times_U F \dar[d] \dar[r] & F \dar[d] \\
    R \Times _U G \dar[r] & G. } $$
  Pour tout schéma $W$,\ifmmode\,\fi{}  on a donc un $2$-morphisme de groupoïdes
  $$ \english\xymatrix @C=4cm {
    \relax\hom(W,\ifmmode\,\fi{} T) \Times_{\hom(W,\ifmmode\,\fi{} F)} \hom(W,\ifmmode\,\fi{} T)
    \ar[d] \dar[r] & 
    \hom(W,\ifmmode\,\fi{} T) \ar[d] 
    \ar[ddl]
    |!{[ld];[d]}\bighole
    _(.3){\hom(W,\ifmmode\,\fi{} \gamma)} \\ 
    \hom(W,\ifmmode\,\fi{} R) \Times_{\hom(W,\ifmmode\,\fi{} U)} \hom(W,\ifmmode\,\fi{} F) 
    \dar[d]
    \dar[r] & 
    \hom(W,\ifmmode\,\fi{} F) 
    \ar@<+2pt>[d]^{\hom(W,\psi)}
    \ar@<-2pt>[d]_{\hom(W,\phi)}
    \\
    \hom(W,\ifmmode\,\fi{} R) \Times _{\hom(W,\ifmmode\,\fi{} U)} \hom(W,\ifmmode\,\fi{} G) 
    \dar[r]
    & \hom(W,\ifmmode\,\fi{} G). } $$
  D'après le lemme \ref{lemme:morobjequiv},\ifmmode\,\fi{}  on a donc 
  $\hom(W,\ifmmode\,\fi{}  \phi) = \hom(W,\ifmmode\,\fi{} \psi)$ pour tout $W$,\ifmmode\,\fi{}  donc 
  $\phi = \psi$.
  
  (d) Montrons que l'application 
  $$ \hom_{ R \dlra U }( F,\ifmmode\,\fi{}  G ) \lra \cat{hom}_{[U/R]} 
  \bigl( [F/R],\ifmmode\,\fi{}  [G/R] \bigr) / \text{isomorphismes}$$
  est surjective. 
  Soit donc 
  $$ \xymatrix{ F/R \ar[rr]^\phi \ar[rd]|*{}="A"
    \enlargexyentry A
    && G/R \ar[ld] 
    \ar@{=>}"A"_\alpha
    \\
    & U/R } $$
  un morphisme de champs au dessus de $U/R$. D'après le lemme
  \ref{lemme:prescarcart},\ifmmode\,\fi{}  
  on a un carré $2$-cartésien 
  $$ \xymatrix{ 
    F \cartesien 
    \ar[r] \ar[d] & F/R \ar[d] \\
    G \ar[r] & G/R.} $$
  Comme dans le diagramme 
  $$\english  \xymatrix{
    F \Times_{F/R} F \ar@{.>}[d] \dar[r] & F \ar[d]_\psi \ar[r]
    \cartesien & F/R \ar[d]^\phi \\
    G \Times_{G/R} G \dar[r] & G \ar[r] & G/R,\ifmmode\,\fi{}  } $$
  il existe une flèche en pointillés (qui fasse tout commuter),\ifmmode\,\fi{}  à
  savoir 
  $$ \xymatrix{ 
    F \Times_{F/R} F \ar[r] & 
    F \Times_{G/R} F \ar[r]^{\psi \times \psi} & 
    G \Times_{G/R} G,\ifmmode\,\fi{}  } $$
  le morphisme 
  $\psi : F \lra G$ est bien équivariant
  et son image est bien (isomorphe au) morphisme $(\phi,\ifmmode\,\fi{} \alpha)$ dont on est
  parti. 
\end{dem}

\begin{rem}\label{rem:fidele}
  \cbstart
  Soient $\cat U = (R \dlra U)$ un groupoïde dans la catégorie des
  ensembles et 
  $(F,\ifmmode\,\fi{}  \pi: F \lra U,\ifmmode\,\fi{}  
  \alpha: R \times_U F \lra F)$ 
  un ensemble $\cat U$-équivariant. 
  $$
  \xymatrix{
    R \times_U F \ar[r]^-\alpha & F
    \ar[d]^\pi 
    \\
    R \ar@<+2pt>[r]^s \ar@<-2pt>[r]_b 
    & U }
  $$
  Cet objet équivariant définit un groupoïde 
  $\cat F = (R \times _U F \dlra F)$,\ifmmode\,\fi{}  dont les morphismes source,\ifmmode\,\fi{}  but,\ifmmode\,\fi{} 
  unité,\ifmmode\,\fi{}  inverse et multiplication sont respectivement 
  $\pr_2$,\ifmmode\,\fi{}  $\alpha$,\ifmmode\,\fi{}  
  $$ 
  e : \definefunction{F}{R\times_U F}{f}{(\id_{\pi(f)},\ifmmode\,\fi{}  f)} 
  \qquad
  i : \definefunction{R\times_U F}{R\times_U F}%
  {(\phi,\ifmmode\,\fi{} f)}{\bigl(\phi^{-1},\ifmmode\,\fi{}  \alpha(\phi,\ifmmode\,\fi{} f)\bigr)}
  $$
  $$ 
  m : 
  \left\{
    \raisebox{.5\depth}{\xymatrix{
        (R \TIMES{s}{U}{\pi} F) 
        \TIMES{\pr_2}{F}{\alpha}
        (R \TIMES{s}{U}{\pi} F) 
        \eq[r] 
        &
        R \TIMES sUb 
        R \TIMES sU{\pi} F
        \ar[r] &
        R \TIMES sU{\pi} F
        \\
        \bigl( \psi,\ifmmode\,\fi{}  \alpha(\phi,\ifmmode\,\fi{} f),\ifmmode\,\fi{}  \phi,\ifmmode\,\fi{}  f \bigr)
        \ar@{|->}[r] &
        (\psi,\ifmmode\,\fi{} \phi,\ifmmode\,\fi{} f) 
        \ar@{|->}[r] &
        (\psi\circ \phi,\ifmmode\,\fi{}  f).
        }}
  \right.
  $$

  On remarquera que $\cat F \lra \cat U$ est un revêtement de
  groupoïdes,\ifmmode\,\fi{}  \ie,\ifmmode\,\fi{}  pour tout objet 
  $f \in \ob \cat F$,\ifmmode\,\fi{}  on a une bijection 
  $$ 
  \definefunction
  {\{\,\ifmmode\,\fi{} \text{flèches de $\cat F$ partant de }f\,\ifmmode\,\fi{} \}}%
  {\{\,\ifmmode\,\fi{} \text{flèches de $\cat G$ partant de }\pi(f)\,\ifmmode\,\fi{} \}}%
  {\xymatrix{
      f \ar[r]^-{(\phi,\ifmmode\,\fi{} f)} & \alpha(\phi,\ifmmode\,\fi{} f)
      }}%
  {\xymatrix{
      \pi f \ar[r]^-\phi & \pi \alpha(\phi,\ifmmode\,\fi{} f).
      }}%
  $$
  En particulier,\ifmmode\,\fi{}  le foncteur $\cat F \lra \cat U$ est fidèle. 

  On remarquera aussi qu'un morphisme d'ensembles équivariants induit
  un morphisme de groupoïdes. 
  \cbend
\end{rem}

\begin{lemme}\label{lemme:morobjequiv}%
  \cbstart
  Soient 
  $\cat U = (R \dlra U)$ un groupoïde (dans la catégorie des
  ensembles),\ifmmode\,\fi{}  $F$ et $G$ des ensembles 
  $\cat U$-équivariants,\ifmmode\,\fi{}  qui définissent des groupoïdes 
  $\cat F = (R \times_U F \dlra F)$
  et 
  $\cat G = (R \times_U G \dlra G)$
  comme dans la remarque précédente. 
  Soient $T$ un ensemble et $h : T \lra F$ une application,\ifmmode\,\fi{}  qui définit
  un groupoïde 
  $\cat T = ( T \times_F T \dlra T)$ et un morphisme de groupoïdes,\ifmmode\,\fi{} 
  encore noté $h$,\ifmmode\,\fi{}  
  $\cat T \lra \cat F$. 
  Soient enfin $\phi$ et $\psi : F \lra G$ deux applications
  équivariantes. 
  S'il existe un 2-morphisme de groupoïdes 
  $\gamma : \phi h \lra \psi h$
  $$ 
  \english
  \xymatrix @!0 @R=1.5pc @C=4pc {
    & \cat F \ar[rd]^\phi \ar@{=>}[dd]^\gamma \\
    \cat T \ar[ru]^h \ar[rd]_h 
    && \cat G \\
    & \cat F \ar[ru]_\psi
    }
  $$
  tel que le diagramme 
  $$ 
  \english
  
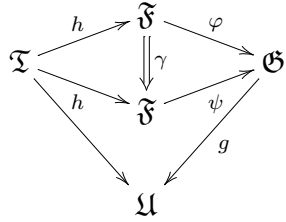
\xymatrix @!0 @R=1.5pc @C=4pc {
    & \cat F \ar[rd]^\phi \ar@{=>}[dd]^\gamma \\
    \cat T \ar[ru]^h \ar[rd]_h 
    \ar[rddd] 
    && \cat G 
    \ar[lddd]^g \\
    & \cat F \ar[ru]_\psi \\ \\
    & \cat U
    }
  $$
  soit 2-commutatif,\ifmmode\,\fi{}  \ie,\ifmmode\,\fi{}  
  $g \gamma = \id_{g \phi h} = \id_{g \psi h}$,\ifmmode\,\fi{}  
  alors 
  $\phi = \psi$ et $\gamma = \id$. 
  \cbend
\end{lemme}
\begin{dem}
  La situation est la suivante. 
  $$ \english
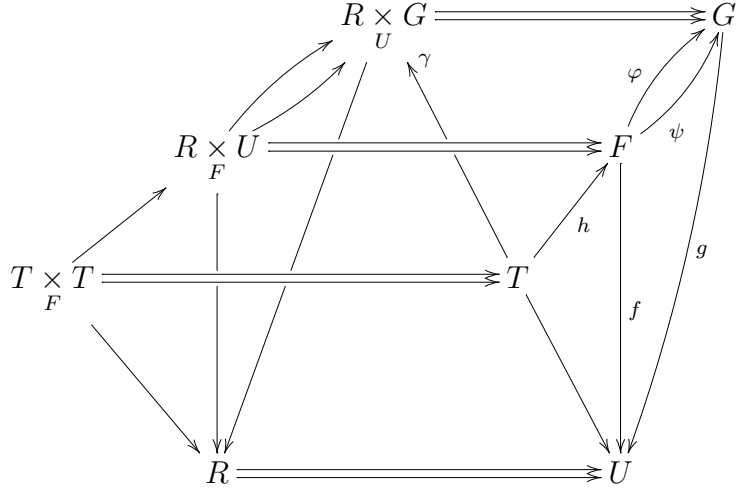
\xymatrix{
    && R \Times_U G \dar[rrr] 
    \ar[ldddd] 
    |!{[ld];[rrd]}\hole
    |!{[lldd];[rdd]}\hole
    &&& G \ar@/^/[ldddd]^g \\
    & R \Times_F U \ar[ddd]
    |!{[ld];[rrd]}\hole
    ^(.3){} \dar[rrr] 
    \ar@/^/[ru] \ar@/_/[ru] &&& 
    F \ar@/^/[ru]^(.4)\phi \ar@/_/[ru]_(.3)\psi \ar[ddd]^f \\
    T \Times_F T \ar[rdd] \dar[rrr] \ar[ru] &&& 
    T \ar[luu]
    |!{[llu];[ur]}\hole
    _(.8){\gamma} \ar[ru]_h \ar[rdd] \\
    \\
    & R \dar[rrr] &&& U
    } $$
  On regarde $T$,\ifmmode\,\fi{}  $G$ et $F$ comme des catégories. 
  Soit $x \in T$. Considérons le foncteur $g$ et plus particulièrement
  l'application 
  $$ g : \definefunction
  {\hom_G( \phi h x,\ifmmode\,\fi{}  \psi h x)}%
  {\hom_U( f h x,\ifmmode\,\fi{}  f h x)}%
  {(\alpha : \phi h x \lra \psi h x)}%
  {(g \alpha : f h x \lra f h x).}%
  $$
  En particulier,\ifmmode\,\fi{}  l'image de 
  $$ \gamma _x : \phi h x \lra \psi h x $$
  est 
  $$ g \gamma _x : f h x \lra f h x. $$
  Mais comme $g \gamma = \id$,\ifmmode\,\fi{}  on a 
  $g \gamma _x = \id_{ f h x }$.
  \cbstart
  D'autre part,\ifmmode\,\fi{}  
  comme le foncteur $g$ est fidèle d'après la remarque
  \ref{rem:fidele},\ifmmode\,\fi{}  on a 
  $\gamma _x = \id_{\phi h x }$ et 
  $\phi h x = \psi h x$.
  \cbend
\end{dem}

\begin{rem}
  \cbstart
  En particulier,\ifmmode\,\fi{}  pour $T=F$ et $h = \id$,\ifmmode\,\fi{}  le lemme nous dit que les
  seuls $2$-morphismes de groupoïdes $\cat F \lra \cat G$ au dessus de
  $\cat U$ sont les identités.
  \cbend
\end{rem}

\begin{lemme}\label{lemme:prescarcart}%
  Soient $X$ un champ algébrique,\ifmmode\,\fi{}  
  $R \dlra U$ une de ses présentations
  et 
  \cbstart
  $(F,\ifmmode\,\fi{}  \pi: F \lra U,\ifmmode\,\fi{}  \alpha : R \times_U F \lra F)$ 
  \cbend
  un revêtement équivariant de $R \dlra U$. 
  On a alors un carré $2$-cartésien
  $$ \xymatrix{
    F \ar[r] \ar[d]_{\pi} \cartesien & \relax[F/R] \ar[d] \\
    U \ar[r] & X. } $$  
\end{lemme}
\begin{dem}
  (a) Commençons par démontrer que l'on a un tel carré $2$-cartésien
  dans la $2$-catégorie des groupoïdes
  si $R \dlra U$ est un groupoïde et $F$ un ensemble équivariant.
  Le $2$-produit-fibré 
  $(F/R) \dtimes_{U/R} U$ a pour objets 
  $$ \biggl\{\,\ifmmode\,\fi{}  
  (f,\ifmmode\,\fi{} u,\ifmmode\,\fi{} \gamma) \in F \times U \times R \ :\ \gamma \in \hom\bigl( u,\ifmmode\,\fi{}  \pi(f)
  \bigr) \,\ifmmode\,\fi{} \biggr\} $$
  et pour morphismes 
  \begin{align*}
    \hspace{2cm}&\hspace{-2cm}
    \hom\bigl( (f,\ifmmode\,\fi{} u,\ifmmode\,\fi{} \gamma),\ifmmode\,\fi{}  (f',\ifmmode\,\fi{} u',\ifmmode\,\fi{} \gamma') \bigr) 
    \\
    &= \left\{\,\ifmmode\,\fi{}  (\phi,\ifmmode\,\fi{} \psi) \in \hom(u,\ifmmode\,\fi{} u') \times \hom(f,\ifmmode\,\fi{} f') 
    \ :\ 
    \raisebox{.5\depth}{\xymatrix{
        u \ar[r]^{\gamma} \ar[d]_\phi & 
        \pi(f) \ar[d]^{\pi(\psi)} \\
        u' \ar[r]_{\gamma'} & 
        \pi(f') }}
    \,\ifmmode\,\fi{} \right\} 
    \displaybreak[0]\\
    &= \bigl\{\,\ifmmode\,\fi{}  (\phi,\ifmmode\,\fi{} \psi) \in \hom(u,\ifmmode\,\fi{} u') \times \hom(f,\ifmmode\,\fi{} f') 
    \ :\ 
    \phi = \id,\ifmmode\,\fi{} 
    \\
    &\hspace{3cm}
    \psi = (\chi,\ifmmode\,\fi{}  f) \in R \Times_U F,\ifmmode\,\fi{} 
    \ \ 
    \alpha(\chi,\ifmmode\,\fi{} f) = f' 
    \text{\ \ et\ \ }
    \chi \gamma = \gamma'
    \,\ifmmode\,\fi{} \bigr\} 
    \displaybreak[0]    \\
    &= \left\{\,\ifmmode\,\fi{}  \chi \in R \ :\ 
      u = u',\ifmmode\,\fi{}  \quad
      s(\chi) = f,\ifmmode\,\fi{}  \quad 
      \alpha(\chi,\ifmmode\,\fi{} f) = f'
      \text{\ \ et\ \ } \chi = \gamma' \gamma^{-1}
      \,\ifmmode\,\fi{} \right\}
    \\
    &=
    \begin{cases}
      * &\text{si } 
      u = u',\ifmmode\,\fi{}  \quad \alpha( \gamma' \gamma^{-1},\ifmmode\,\fi{}  f) = f' \\
      & \hspace{1cm} 
      \text{ (et }
      s(\gamma'\gamma^{-1}) = \pi(f)
      \text{,\ifmmode\,\fi{}  mais c'est automatique)}
      \\
      \nothing &\text{sinon}
    \end{cases}
  \end{align*}
  On constate qu'il s'agit d'un ensemble muni d'une relation
  d'équivalence :
  $$ (f,\ifmmode\,\fi{}  u,\ifmmode\,\fi{}  \gamma) \sim (f',\ifmmode\,\fi{} u',\ifmmode\,\fi{} \gamma') 
  \ssi \alpha( \gamma' \gamma^{-1},\ifmmode\,\fi{}  f) = f'. $$
  On remarque que tout élément $(f,\ifmmode\,\fi{} u,\ifmmode\,\fi{} \gamma)$ est équivalent à un
  élément de la forme 
  $(f',\ifmmode\,\fi{}  \pi(f'),\ifmmode\,\fi{}  \id)$ :
  $$ (f,\ifmmode\,\fi{}  u,\ifmmode\,\fi{}  \gamma) \sim \bigl( \alpha(\gamma^{-1},\ifmmode\,\fi{} f),\ifmmode\,\fi{}  \pi(
  \alpha(\gamma^{-1},\ifmmode\,\fi{} f)),\ifmmode\,\fi{}  \id\bigr). $$
  D'autre part,\ifmmode\,\fi{}  deux élément de la forme 
  $(f,\ifmmode\,\fi{} \pi(f),\ifmmode\,\fi{} \id)$ sont équivalents \ssi ils sont égaux. 
  On a donc 
  $$ F/R \Dtimes_{U/R} U \equiv F.$$

  (b) On peut supposer $X = [U/R]$. Soit
  maintenant $T$ un schéma.
  Nous allons montrer que le carré 
  $$ \xymatrix{
    F_T \ar[r] \ar[d] \cartesien & \relax[F/R]_T \ar[d] \\
    U_T \ar[r] & X_T } $$    
  est $2$-cartésien en exhibant une équivalence de catégories 
  $$ F_T \equiv U_T \Times_{[U/R]_T} [F/R]_T. $$
  Mais d'après \cite{vistoli}, 
  \begin{align*}
    \relax[U/R]_T &= \dlimind\limits_{T' \lra T} 
    [U/R]_{T' \lra T} \\
    \text{où }\relax[U/R]_{T' \lra T}
    &=\cat{hom}
    \left(
      \raisebox{.5\depth}{\xymatrix{
          T' \Times_T T' \dar[d] 
          \ar@{}[rd]|{\displaystyle ,} &
          R \dar[d]
          \\ 
          T' & U
          }}
    \right).
  \end{align*}
  Nous allons donc simplement montrer que l'on a un carré $2$-cartésien de
  groupoïdes 
  $$ \xymatrix{
    F_{T'\lra T} \ar[r] \ar[d] \cartesien & \relax[F/R]_{T'\lra T} \ar[d] \\
    U_{T'\lra T} \ar[r] & X_{T'\lra T} } $$
  où $F_{T'\lra T}$ et $U_{T'\lra T}$ sont des ensembles

  Pour cela,\ifmmode\,\fi{}  il suffit de montrer que le groupoïde
  $$\english
  G_1 = U_{T' \lra T} / R_{T' \lra T} 
  \,\ifmmode\,\fi{}  : \quad 
  \cat{hom}
  \left(
    \raisebox{.5\depth}{\xymatrix{
        T' \Times_T T' \dar[d] 
        \ar@{}[rd]|{\displaystyle ,} &
        R \dar[d]
        \\ 
        T' & R
        }}
  \right)
  \dlra 
  \cat{hom}
  \left(
    \raisebox{.5\depth}{\xymatrix{
        T' \Times_T T' \dar[d] 
        \ar@{}[rd]|{\displaystyle ,} &
        U \dar[d]
        \\ 
        T' & U
        }}
  \right)
  $$
  qui apparait dans (a) est équivalent à
  $$
  \english
  G_2 = X_{T' \lra T} = \,\ifmmode\,\fi{}  \cat{hom}
  \left(
    \raisebox{.5\depth}{\xymatrix{
        T' \Times_T T' \dar[d] 
        \ar@{}[rd]|{\displaystyle ,} &
        R \dar[d]
        \\ 
        T' & U
        }}
  \right)
  $$ et d'appliquer (a). 
  
  À un objet du premier groupoïde,\ifmmode\,\fi{}  $G_1$,\ifmmode\,\fi{}  \ie,\ifmmode\,\fi{}  à
  un morphisme de groupoïdes
  $$
  \xymatrix{T' \Times_T T' \ar[d] \dar[r] & T' \ar[d] \\
    U \dar[r] & U} $$
  on associe un objet du
  second,\ifmmode\,\fi{}  $G_2$,\ifmmode\,\fi{}  par composition
  $$
  \xymatrix{T' \Times_T T' \ar[d] \dar[r] & T' \ar[d] \\
    U \ar[d] \dar[r] & U\ar[d] \\
    R \dar[r] & U.} $$
  Cette application (entre
  les ensembles d'objets de nos deux groupoïdes)
  est bijective car les seuls morphismes du
  groupoïde
  $$
  T' \Times_T T' \dlra T'$$
  sont les identités
  et un morphisme de groupoïdes
  $$
  \xymatrix{T' \Times_T T' \ar[d] \dar[r] & T' \ar[d] \\
    R \dar[r] & U} $$
  préserve les identités et se
  factorise donc par
  $$
  \xymatrix{T' \Times_T T' \ar[d] \dar[r] & T' \ar[d] \\
    U \ar[d] \dar[r] & U\ar[d] \\
    R \dar[r] & U.} $$
  
  Regardons maintenant les morphismes de nos
  groupoïdes.
  \begin{align*}
    &\hspace{-2cm}
    \hom_{G_1} \left(
      \raisebox{.5\depth}{\xymatrix{
          T' \Times_T T' \ar[d] \dar[r] & T' \ar[d]^f 
          \ar@{}[rd]|{\displaystyle\text{,}}
          & 
          T' \Times_T T' \ar[d] \dar[r] & T' \ar[d]^g
          \\
          U \dar[r] & U 
          &
          U \dar[r] & U 
          }}
    \right)
    \\
    &= \left\{
      \raisebox{.5\depth}{\xymatrix{
          T' \Times_T T' \ar[d] \dar[r] & T' \ar[d]^\alpha
          \\ R \dar[r] & R }}
      \text{ tels que }
      s \alpha = f 
      \text{ et }
      b \alpha = g
    \right\} 
    \displaybreak[0]\\
    &= \left\{
      \raisebox{.5\depth}{\xymatrix{
          T' \Times_T T' \ar[d] \dar[r] & T' \ar[d]^\alpha
          \\ R \ar[d] \dar[r] & R 
          \ar@<+2pt>[d]^b \ar@<-2pt>[d]_s 
          \\ R \dar[r] & U }}
      \text{ tels que }
      s \alpha = f 
      \text{ et }
      b \alpha = g
    \right\}
    \displaybreak[0]\\
    &= \left\{
      \raisebox{.5\depth}{\xymatrix{
          T' \Times_T T' \ar[d] \dar[r] & T' \ar[ld]_\alpha
          \ar@<+2pt>[d]^g \ar@<-2pt>[d]_f
          \\ R \ar@<+2pt>[r]^s \ar@<-2pt>[r]_b & U
          }}
      \text{ tels que }
      s \alpha = f 
      \text{ et }
      b \alpha = g
    \right\}
    \displaybreak[0]\\    
    &= \hom_{G_2}\left(
      \raisebox{.5\depth}{\xymatrix{
          T' \Times_T T' \ar[d] \dar[r] & T' \ar[d]^f 
          \ar@{}[rd]|{\displaystyle\text{,}} &  
          T' \Times_T T' \ar[d] \dar[r] & T' \ar[d]^g
          \\
          R \dar[r] & U 
          &
          R \dar[r] & U 
          }}
      \right)  
  \end{align*}
\end{dem}

\section{Théorème de Van Kampen}


\def\DAR{
  \renewcommand{\dar}[3][r]{
    \ar@<+2pt>[##1]^{##2}
    \ar@<-2pt>[##1]_{##3}
    }
}

Nous rappelons le théorème de Van Kampen tel qu'il est présenté dans 
\cite{giraud} ou \cite[IV.5]{SGA1}, avant de le traduire dans un
langage plus concret en termes de carré cocartésien de progroupoïdes. 

\begin{defi}
  Soit $\f F : \cat F \lra \cat C$ une catégorie fibrée et 
  $\phi : X \lra Y$ un morphisme de $\cat C$. 
  On note $p_i$ et $p_{ij}$ les différentes projections :
  $$\english
  \xymatrix{
    *!!<0pt,.7ex>+!R!<-5pt,0pt>{X\Times_Y X \Times_Y X}
    \ar@/^1pc/[r]^-{p_{12}}
    \ar[r]    ^-{p_{13}}
    \ar@/_1pc/[r]_-{p_{23}}
    & \relax\smash[b]{X \Times_Y X} 
    \ar@<+2pt>[r]^-{p_1}
    \ar@<-2pt>[r]_-{p_2}
    & X \ar[r]^\phi 
    & Y.
    }
  $$
  La \d{catégorie des données de descente} pour $\f F$ relatives à
  $\phi$,\ifmmode\,\fi{}  
  notée $\desc\phi$,\ifmmode\,\fi{} 
  a pour objets les 
  $(F,\ifmmode\,\fi{}  \alpha)$,\ifmmode\,\fi{}  où 
  $F \in \ob \cat F_X$ est un objet de $\cat F$ au dessus de $X$ et 
  $\alpha : p_1^* F \lra p_2^* F$ est un isomorphisme au dessus de 
  $\id_{X\times_YX}$ tel que 
  $$
  \english
  \renewcommand{\eq}[2]
  {\ar@<-3pt>@{-}"#1";"#2"
    \ar@<-1pt>@{}"#1";"#2"|<{}="gauche"
    \ar@<+0pt>@{}"#1";"#2"|-{}="milieu"
    \ar@<+1pt>@{}"#1";"#2"|>{}="droite"
    \ar@/^2pt/@{-}"gauche";"milieu"
    \ar@/_2pt/@{-}"milieu";"droite"
    \ar@{}"#1";"#2"}
  \def\MyNode{\ifcase\xypolynode\or
     p_{23}^* p_2^* F
   \or
     p_{23}^* p_1^* F
   \or
     p_{12}^* p_2^* F
   \or
     p_{12}^* p_1^* F
   \or
     p_{13}^* p_1^* F
   \or
     p_{13}^* p_2^* F
   \fi
 }%
 \xy/r5pc/:
   \xypolygon6{~>{}\txt{\ \ \strut\ensuremath{\MyNode}}} 
   \ar "2";"1" ^{p_{23}^*\alpha}
   \ar "4";"3" ^{p_{12}^*\alpha}
   \ar "5";"6" _{p_{13}^*\alpha}
   \eq32 \eq45 \eq61
 \endxy
 $$
 et dont les morphismes 
 $(F,\ifmmode\,\fi{}  \alpha) \lra (F',\ifmmode\,\fi{}  \alpha')$ sont les morphismes 
 $f : F \lra F'$ dans $\cat F_X$ tels que 
 $$\xymatrix{ p_1^* F \ar[r]^\alpha \ar[d]_{p_1^* f} &
   p_2^* F \ar[d] ^{p_2^* f} \\
   p_1^* F' \ar[r]_{\alpha'} & p_2^* F'.}$$
\end{defi}

\begin{rem}
  Soit $F : \cat F \lra \cat C$ une catégorie fibrée et 
  $\phi : X \lra Y$ un morphisme de $\cat C$. 
  D'après \cite[IV.5]{SGA1}, on a 
  une équivalence de
  catégories 
  $$ \desc \phi \equiv \dlimpro \left(
    \cat F _X 
    \dlra \cat F_{X \Times_Y X}
    \tlra \cat F_{X \Times_Y X \Times_Y X}
  \right) $$
  où le $2$-système projectif est défini par :
  \begin{theoenum}
  \item trois catégories $\F1$,\ifmmode\,\fi{}  $\F2$ et $\F3$ ;
  \item des foncteurs entre ces catégories 
    $$\english 
    \xymatrix{\relax
      \F1
      \ar@<+2pt>[r]^-{p^*_1}
      \ar@<-2pt>[r]_-{p^*_2}
      & \F2 
      \ar@<+4pt>[r]^-{p^*_{ij}}
      \ar[r]
      \ar@<-4pt>[r]
      & \F3 ;} $$
  \item des $2$-isomorphismes 
    $$\english\xymatrix @!0 @R=3pc @C=4pc {\relax
      & \F2 \ar[rd]^{p^*_{12}} \ar@{=>}[dd]^{\theta_1} \\
      \F1 \ar[ru]^{p^*_1} \ar[rd]_{p^*_2} && \F3 \\
      & \F2 \ar[ru]_{p^*_{13}} } 
    \qquad
    \english\xymatrix @!0 @R=3pc @C=4pc {\relax
      & \F2 \ar[rd]^{p^*_{12}} \ar@{=>}[dd]^{\theta_2} \\
      \F1 \ar[ru]^{p^*_1} \ar[rd]_{p^*_2} && \F3 \\
      & \F2 \ar[ru]_{p^*_{23}} } 
    $$
    $$
    \english\xymatrix @!0 @R=3pc @C=4pc {\relax
      & \F2 \ar[rd]^{p^*_{13}} \ar@{=>}[dd]^{\theta_3} \\
      \F1 \ar[ru]^{p^*_1} \ar[rd]_{p^*_2} && \F3 \\
      & \F2. \ar[ru]_{p^*_{23}} } 
    $$      
  \end{theoenum}
\end{rem}

\begin{defi}
  Soit $\f F : \cat F \lra \cat C$ une catégorie
  fibrée. On dit qu'un morphisme $\phi : X \lra Y$
  de $\cat C$ est un \d{morphisme de descente
    effective} si le foncteur
  $$\definefunction{\cat F _X}{\desc \phi}%
  {F}{ ( \phi^* F,\ifmmode\,\fi{}  \ 
    p_1^* \phi^* F \lra (p_1 \phi)^* F = (p_2 \phi)^* F \lra p_2^*
    \phi^* F ) }$$
  est une équivalence de catégories.
\end{defi}

\begin{theo}\label{theo:ddgir}%
  Un morphisme couvrant 
  $\phi : X \lra *$ dans un topos $\T$ est un morphisme de descente
  effective pour la catégorie fibrée associée au pseudo-foncteur 
  $$\definefunctor{\T}{\dcat{Cat}}{T}{\T/T}{f}{f^*.}$$
\end{theo}
\begin{dem}
  Voir \cite{giraud}.
\end{dem}

\begin{coro}
  Un morphisme couvrant 
  $\phi : X \lra *$ dans un topos $\T$ est un morphisme de descente
  effective pour la catégorie fibrée associée au pseudo-foncteur 
  $$
  \definefunctor{\T}{\dcat{Cat}}{T}{\SLC\T/T}{f}{f^*}
  \qquad
  \cbstart
  \text{respectivement}
  \qquad
  \definefunctor{\T}{\dcat{Cat}}{T}{\SLCF\T/T}{f}{f^*.}
  \cbend
  $$
\end{coro}
\begin{dem}
  Il suffit de remarquer que 
  \begin{align*}
    (\SLC \T) /X &\iso \SLC (\T / X) \\
    &\iso \SLC \cat{Desc}(\phi,\ifmmode\,\fi{}  \T/\cdot) \\
    &\iso \cat{Desc}(\phi,\ifmmode\,\fi{}  \SLC \T/ \cdot)
  \end{align*}
\cbstart
  et de procéder de même pour le cas de $\SLCF$.
\cbend
\end{dem}

\begin{prop}
  Soit $\T$ un topos,\ifmmode\,\fi{}  
  $U$ et $V$ des objets de $\T$ tels que 
  $ \phi : U \amalg V \lra *$ soit couvrant et que 
  $U \times U \iso U$ et $V \times V \iso V$. On a alors un carré
  $2$-cartésien de catégories et de foncteurs « image réciproque » de
  topos,\ifmmode\,\fi{}  
  $$ \xymatrix{\relax
    \T \ar[r] \ar[d] \cartesien & \T/U \ar[d] \\ \T/V \ar[r] & 
    \T / U \times V.} $$
\end{prop}

\begin{dem}
  On sait déjà que 
  $\T \equiv \desc\phi$ : nous allons expliciter la catégorie
  $\desc\phi$ des données de descente et voir que cette description
  est la même que celle du $2$-produit fibré 
  $\T/U \times_{\T/U\times V} \T/V$ : on pourra alors conclure grace
  au théorème \ref{theo:ddgir}.
  Une donnée de descente relative à $f$ est la donnée de 
  \begin{itemize}
  \item Un objet de 
    $\T/U\amalg V$,\ifmmode\,\fi{}  que l'on peut écrire 
    $F \amalg G$,\ifmmode\,\fi{}  où $F$ et $G$ sont des objets de 
    $\T/U$ et $\T/V$ respectivement ;
  \item Un isomorphisme 
    $ \theta : p_1^* ( F \amalg G) \lra p_2^* ( F \amalg G)$ 
    tel que 
    $$\english \xymatrix @!0 @C=4pc @R=4pc {
      & q_2^* ( F \amalg G ) \ar[rd]^{p_{23}^*\theta} \\
      q_1^* ( F \amalg G )  
      \ar[ru]^{p_{12}^*\theta} 
      \ar[rr]_{p_{13}^*\theta} 
      && q_3^* ( F \amalg G ).} $$
  \end{itemize}
  Or,\ifmmode\,\fi{}  
  on peut écrire $\theta$ sous la forme 
  $$ F|_{U\times U} \amalg F|_{U \times V} \amalg G|_{U \times V}
  \amalg G|_{V\times V}
  \lra
  F \amalg G|_{U \times V} \amalg F|_{U \times V} \amalg G,\ifmmode\,\fi{}  $$
  \ie,\ifmmode\,\fi{}  puisque $F|_{U\times U} = F|_U = F$ et 
  $G|_{V\times V} = G|_V = G$,\ifmmode\,\fi{} 
  $$ F \amalg F|_{U \times V} \amalg G|_{U \times V} \amalg G
  \lra
  F \amalg G|_{U \times V} \amalg F|_{U \times V} \amalg G. $$
  Une donnée de descente est donc la donnée d'isomorphismes
  $$\xymatrix{
    F \ar[r]^\alpha & F,\ifmmode\,\fi{}  &
    F|_{U \times V} \ar[r]^\beta & G|_{U \times V} \\
    G \ar[r]^\delta & G,\ifmmode\,\fi{}  &
    G|_{U \times V} \ar[r]^\gamma & F|_{U \times V}\rlap{{};{}}}$$
  tels que 
  $$ \english
  \xymatrix @!0 @R=3pc @C=3pc {
    & F \ar[rd]^\alpha \\ F \ar[ru]^\alpha \ar[rr]_\alpha && F }
  \qquad
  \xymatrix @!0 @R=3pc @C=3pc {
    & G \ar[rd]^\delta \\ G \ar[ru]^\delta \ar[rr]_\delta && G }
  \qquad
  \english
  \xymatrix @!0 @R=3pc @C=3pc {
    &   G|_{U \times V} \ar[rd]^\gamma \\ 
    F|_{U \times V} \ar[ru]^\beta \ar[rr]_\alpha && 
    F|_{U \times V,\ifmmode\,\fi{} \,\ifmmode\,\fi{} } } 
  $$
  \ie,\ifmmode\,\fi{}  
  $\alpha = \id$,\ifmmode\,\fi{}  
  $\delta = \id$,\ifmmode\,\fi{}  
  $\gamma = \beta^{-1}$.
  On expliciterait de même les morphismes de données de descente. 
  En définitive,\ifmmode\,\fi{}  la catégorie des données de descente 
  est équivalente à la catégorie des 
  $(F,\ifmmode\,\fi{} G,\ifmmode\,\fi{} \beta)$,\ifmmode\,\fi{}  
  où $F$ et $G$ sont des objets de 
  $\T/U$ et $\T/V$ et 
  $\beta$ est un isomorphisme 
  $ F|_{U \times V} \lra G|_{U \times V} $ :
  c'est exactement la description du $2$-produit fibré 
  $\T_U \Times_{\T/U\times V} \T/V$.
\end{dem}

\begin{coro}
  Soit $\T$ un topos,\ifmmode\,\fi{}  
  $U$ et $V$ des objets de $\T$ tels que 
  $ \phi : U \amalg V \lra *$ soit couvrant et que 
  $U \times U \iso U$ et $V \times V \iso V$. On a alors un carré
  $2$-cartésien de catégories et de morphismes image réciproque de
  topos,\ifmmode\,\fi{}  
  $$ \xymatrix{\relax
    \SLC \T \ar[r] \ar[d] \cartesien & \SLC \T/U \ar[d] \\ \SLC \T/V \ar[r] & 
    \SLC \T / U \times V.} $$
\end{coro}

\begin{dem}
  Comme 
  $\SLC(\T/X) = (\SLC \T)/X$,\ifmmode\,\fi{}  il suffit d'apliquer le résultat
  précédent à $\SLC\T$. 
\end{dem}

\begin{coro}
  Soit $\T$ un topos,\ifmmode\,\fi{}  
  $U$ et $V$ des objets de $\T$ tels que 
  $ \phi : U \amalg V \lra *$ soit couvrant et que 
  $U \times U \iso U$ et $V \times V \iso V$. 
  Soient 
  \begin{itemize}
  \item $P_{U\times V}$ un ensemble de points-bases de 
    $\T/U\times V$,\ifmmode\,\fi{}  au moins un dans chaque composante connexe ;
  \item $P_U$ un ensemble de points-bases de 
    $\T/U$,\ifmmode\,\fi{}  comprenant les images des éléments de 
    $P_{U\times V}$,\ifmmode\,\fi{}  avec au moins un point dans chaque composante
    connexe ; 
  \item $P_V$ un ensemble de points-bases de 
    $\T/V$,\ifmmode\,\fi{}  comprenant les images des éléments de 
    $P_{U\times V}$,\ifmmode\,\fi{}  avec au moins un point dans chaque composante
    connexe ; 
  \item $P$ un ensemble de points-bases de 
    $\T$ comprenant les images des éléments de 
    $P_U$ et $P_V$ avec au moins un point dans
    chaque composante connexe.
  \end{itemize}
  On a alors un carré $2$-cocartesien de progroupoïdes fondamentaux,\ifmmode\,\fi{}  
  $$\xymatrix{\relax
    \pi_1( \T\!,\ifmmode\,\fi{} \,\ifmmode\,\fi{}  P ) \cartesien &
    \pi_1( \T/U,\ifmmode\,\fi{} \,\ifmmode\,\fi{}  P_U) \ar[l] \\
    \pi_1( \T/V,\ifmmode\,\fi{} \,\ifmmode\,\fi{}  P_V) \ar[u] & 
    \pi_1( \T/U\times V,\ifmmode\,\fi{} \,\ifmmode\,\fi{}  P_{U\times V}). \ar[u] \ar[l]} $$
\end{coro}
\begin{dem}
  On a un carré $2$-cartésien 
  $$ \xymatrix{ 
    \SLC \T \ar[r] \ar[d] \cartesien & 
    \SLC \T/U \ar[d] \\
    \SLC \T/V \ar[r] &
    \SLC \T/U\times V,\ifmmode\,\fi{} \,\ifmmode\,\fi{}  } $$
  qui peut s'écrire 
  $$ \xymatrix{\relax
    \B \pi_1 \T \ar[r] \ar[d] \cartesien &
    \B \pi_1 U \ar[d] \\
    \B \pi_1 V \ar[r] &
    \B \pi_1 (U \times V) } $$
  où,\ifmmode\,\fi{}  pour simplifier,\ifmmode\,\fi{}  on a noté 
  \begin{align*}
    \pi_1 \T &= \pi_1 (\T\!,\ifmmode\,\fi{}  P) \\
    \pi_1 U &= \pi_1( \T/U,\ifmmode\,\fi{}  P_U) \\
    \pi_1 V &= \pi_1( \T/V,\ifmmode\,\fi{}  P_V) \\
    \pi_1 (U\times V) &= \pi_1( \T/U \times V,\ifmmode\,\fi{}  P_{U\times V}). \\
  \end{align*}
  D'après les choix de points-bases que l'on a effectués,\ifmmode\,\fi{}  ces
  morphismes proviennent de morphismes de groupoïdes 
  $$\xymatrix{
    \pi_1 \T & \pi_1 U \ar[l] \\ \pi_1 V \ar[u] & \pi_1(U \times V).
    \ar[u] \ar[l] } $$
  Remarquons que,\ifmmode\,\fi{}  pour tout groupe $\pi$,\ifmmode\,\fi{}  on a 
  $ \B \pi \equiv \cat{hom}(\pi,\ifmmode\,\fi{}  \cat{Ens}). $
  On a donc 
  \begin{align*}
    \B \left( \pi_1 U \Dsa_{\pi_1 (U \times V)} \pi_1 V \right) 
    &\equiv \cat{hom}( \pi_1 U \Dsa_{\pi_1 (U \times V)} \pi_1 V,\ifmmode\,\fi{} \,\ifmmode\,\fi{} 
    \cat{Ens} ) \\
    &\equiv 
    \cat{hom}( \pi_1 U ,\ifmmode\,\fi{} \,\ifmmode\,\fi{}  \cat{Ens} )
    \DTimes_{ \cat{hom}\bigl( \pi_1 (U \times V) ,\ifmmode\,\fi{} \,\ifmmode\,\fi{}  \cat{Ens} \bigr) }
    \cat{hom}( \pi_1 V,\ifmmode\,\fi{} \,\ifmmode\,\fi{}  \cat{Ens} ) \\
    &\equiv \cat B \pi_1 U \Dtimes_{\cat B \pi_1 (U \times V)}
    \cat B \pi_1 V \\
    &\equiv \SLC \T \\
    &\equiv \cat B \pi_1 \T.
  \end{align*}
  Donc,\ifmmode\,\fi{}  
  $$ \pi_1(\T\!,\ifmmode\,\fi{} \,\ifmmode\,\fi{}  P) \equiv 
  \pi_1 (\T/U,\ifmmode\,\fi{} \,\ifmmode\,\fi{}  P_U) 
  \Dsa_{\pi_1(\T/U\times V,\ifmmode\,\fi{} \,\ifmmode\,\fi{}  P_{U\times V})}
  \pi_1(\T/V,\ifmmode\,\fi{} \,\ifmmode\,\fi{}  P_V). \qed $$
\end{dem}

Nous allons maintenant voir que dans certains cas,\ifmmode\,\fi{}  cette $2$-somme
amalgamée peut se calculer comme une $1$-somme amalgamée. 
\begin{lemme}
  Soient des morphismes de groupoïdes
  $$ \xymatrix{ C \ar[r]^a \ar[d]_b & A \\ B } $$
  qui établissent des bijections entre les ensembles d'objets. 
  La somme amagamée $A \ast_C B$ est aussi une $2$-somme amalgamée. 
\end{lemme}
\begin{dem}
  Nous allons montrer que $A \ast_C B$ vérifie la propriété
  $2$-universelle d'une $2$-somme amalgamée. 
  Considérons une présentation de chacun de ces groupoïdes
  \begin{align*}
    A &= \bigl\langle (a_i)_{i\in I} ; (r_i)_{i\in I'} \bigr\rangle \\
    B &= \bigl\langle (b_j)_{j\in J} ; (s_i)_{i\in J'} \bigr\rangle \\
    C &= \bigl\langle (c_k)_{k\in K} ; (t_i)_{i\in K'} \bigr\rangle \\
    A\Ast_CB &= \Bigl\langle (a_i)_{i\in I} ,\ifmmode\,\fi{}  (b_j)_{j\in J} ; 
    (r_i)_{i\in I'},\ifmmode\,\fi{}  (s_j)_{j\in J'},\ifmmode\,\fi{}  
    \bigl( a(c_k) b(c_k)^{-1} \bigr)_{k \in K}
    \Bigr\rangle
  \end{align*}
  et un carré $2$-commutatif
  $$ \xymatrix{ C \ar[r]^a \ar[d]_b & A \ar[d]^f \\ B \ar[r]_g
    \ar@{=>}[ru] & D. } $$
  Montrons qu'il existe un diagramme 
  $$ \english \xymatrix{ 
    C \ar[r]^a \ar[d]_b & A \ar[d]_{a'} \ar@/^2pc/[rdd]^f|*{}="B" \\
    B \ar[r]^{b'} \ar@/_2pc/_g[rrd]|*{}="A" 
    \enlargexyentry A
    \enlargexyentry B
    & 
    A \Ast_C B \ar@{.>}[rd]^h \ar@{=>}"B"_-{\beta'} \\ && D. 
    \ar@{=>}"A";[lu]^-{\alpha'}
    } $$
  On peut définir le morphisme $h$ par 
  \begin{align*}
    h(x) &= f(x) &&\hspace{-2cm} \text{si } x \in \ob C = \ob A \\
    h(a_i) &= f(a_i) &&\hspace{-2cm} \text{si } i \in I \\
    h(b_j : x \lra y) &= \alpha_y g(b_j) \alpha_x^{-1} 
    &&\hspace{-2cm} \text{si } j \in J.
  \end{align*}
  C'est bien défini,\ifmmode\,\fi{}  car les relations $h(r_i) = \id$ (pour $i\in I'$) et 
  $h(s_j)=\id$ (pour $j\in J'$) sont vérifiées car elles le sont dans $A$ ou $B$
  et les relations 
  $ h( a(c_k) ) = h( b(c_k) ) $ (pour $k\in K'$) sont vérifiées car
  $\alpha$ est une transformation naturelle : 
  $$ \xymatrix{
    x \ar[d]_{c_k} & g(x) \ar[r]^{\alpha x} \ar[d]_{ gbc_k } &
    fx \ar[d]^{ fac_k } \\ 
    y & gy \ar[r]_{\alpha y} & fy. } $$
  On peut définir $\alpha'$ et $\beta'$ par 
  $\beta' = \id$ et $\alpha'_x = \alpha_x$ (pour $x \in \ob B = \ob C$).
  On obtient bien un diagramme $2$-commutatif,\ifmmode\,\fi{}  \ie,\ifmmode\,\fi{}  
  $$ \xymatrix{ 
    gb \ar[r]^{\alpha' b} \ar[d]_\alpha & 
    h b' b \ar@{=}[d] \\ ha & h a' a. \ar[l]^{\beta' a} } $$
  
  Considérons un autre diagramme $2$-commutatif
  $(h',\ifmmode\,\fi{} \beta,\ifmmode\,\fi{} \gamma)$,\ifmmode\,\fi{} 
  $$ \english \xymatrix{ 
    C \ar[r]^a \ar[d]_b & A \ar[d]_{a'} \ar@/^2pc/[rdd]^f|*{}="B" \\
    B \ar[r]^{b'} \ar@/_2pc/_g[rrd]|*{}="A"
    \enlargexyentry A
    \enlargexyentry B
    &
    A \Ast_C B \ar@{.>}[rd]^{h'} \ar@{=>}"B"^-{\gamma} \\ && D. 
    \ar@{=>}"A";[lu]^-{\beta}
    } $$
  On veut montrer qu'il existe un unique 
  $$ \xymatrix{ \llap{$A \Ast_C{}$} B \UN{h}{h'}{\delta} & D } $$
  tel que 
  $$ \raisebox{.5\depth}{\xymatrix{ 
      h a' \ar[rd]_{\id} \ar[r]^{\delta a'} & h' a'
      \ar[d]^\gamma \\ & f }}
  \qquad\text{et}\qquad
  \raisebox{.5\depth}{\xymatrix{ 
      g \ar[r]^\alpha \ar[rd]_\beta & h b' \ar[d]^{\delta b'}
      \\ & h' b'}}
  \qquad
  \text{commutent.}
  $$
  Si un tel $\delta$ existe,\ifmmode\,\fi{}  on a nécessairement,\ifmmode\,\fi{}  pour tout objet $x$
  de $A$,\ifmmode\,\fi{}  
  $ \gamma_x \circ \delta_{a' x} = \id_x,\ifmmode\,\fi{}  $
  \ie,\ifmmode\,\fi{}  
  $\delta_{a' x} = \gamma_x^{-1}$ ;
  or,\ifmmode\,\fi{}  comme $a'$ est une bijection sur les objets,\ifmmode\,\fi{}  on peut
  identifier~$x$ 
  et~$a' x$ et écrire 
  $$ \delta_x = \gamma_x^{-1}. $$
  Ceci montre que $\delta$ est unique. Vérifions que $\delta$,\ifmmode\,\fi{}  ainsi
  défini,\ifmmode\,\fi{}  est une transformation naturelle : il s'agit de montrer que
  pour tout morphisme 
  $\phi : x \lra y$ de 
  $A \ast_C B$,\ifmmode\,\fi{}  on a 
  $$
  \raisebox{.5\depth}{\xymatrix{
      x \ar[d]_\phi & hx \ar[d]_{h \phi} \ar[r]^{\delta x} & h' x
      \ar[d]^{h' \phi} \\ 
      y & h y \ar[r]_{\delta y} & h' y }}
  \qquad\text{\ie,\ifmmode\,\fi{} }\qquad
  \raisebox{.5\depth}{\xymatrix{
      fx \ar[d]_{h \phi} & h' x \ar[d]^{h' \phi} \ar[l]_{\gamma x} \\ f
      y & h' y. \ar[l]^{\gamma y} }}
  $$
  Il suffit de l'établir pour les générateurs $a_i$ et $b_j$ de $A
  \ast_C B$. 
  D'une part,\ifmmode\,\fi{}  on a bien
  $$ \xymatrix{ 
    fx \ar[d]_{f a_i} & h' x \rlap{${}=h' a' x$} \ar[d]^{h' a_i = h'
      a' a_i} \ar[l]_{\gamma x} \\ fy & 
    h' y \rlap{${}=h' a' y$} \ar[l]^{\gamma y} } $$
  car $\gamma$ est une transformation naturelle. 
  D'autre part,\ifmmode\,\fi{}  comme le diagramme $(h',\ifmmode\,\fi{}  \beta,\ifmmode\,\fi{}  \gamma)$ est
  $2$-commutatif,\ifmmode\,\fi{}  \ie,\ifmmode\,\fi{}  
  $$ \xymatrix{ 
    gbx \ar[r]^{\alpha x} \ar[d]_{\beta b x} & f a x \\ 
    h' b' b x \ar@{=}[r] & h' a' a x, \ar[u]_{\gamma_{ax} } }
  $$
  et comme $\beta$ est une transformation naturelle,\ifmmode\,\fi{}  \ie,\ifmmode\,\fi{}  
  $$\xymatrix{ 
    g b x \ar[r]^{\beta b x} \ar[d] _{g b_j} & h' b' b x \ar[d]^{h'
      b_j} \\ g b y \ar[r]_{\beta b y} & h' b' b y,\ifmmode\,\fi{} } 
  $$
  on a bien 
  $$\xymatrix{ 
    g b x \ar[d]_{g b_j} \ar[r]^{\alpha x} 
    & f a x & h' a' a x \ar[l]_{\gamma a' x} \ar@{=}[r] & 
    h' b' b x \ar[d]^{h' b_j} \\
    g b y \ar[r]_{\alpha y} & f a y & h' a' a y \ar[l]^{\gamma a' y}
    \ar@{=}[r] & h' b' b y. } $$
  Par conséquent,\ifmmode\,\fi{}  $\gamma$ est bien une transformation naturelle. 
  Enfin,\ifmmode\,\fi{}  on a 
  $$ \xymatrix{ h a' \ar[rd]_{\id} \ar[r]^{\delta a'} & h' a'
    \ar[d]^\gamma \\ & f }
  $$
  par définition de $\delta$ et montrer 
  $$
  \xymatrix{ g \ar[r]^\alpha \ar[rd]_\beta & h b' \ar[d]^{\delta b'}
    \\ & h' b'}
  $$
  se ramène à établir
  $$ \xymatrix{ 
    gx \ar[r]^{\alpha x} \ar[rrd]_{\beta x} &
    h b' x \ar@{=}[r] & h a' x \ar@{=}[r] & f x \\
    && h' b' x \ar@{=}[r] & h' a x \ar[u]_{\gamma x} } $$
  qui n'est autre que la $2$-commutativité du diagramme 
  $(h',\ifmmode\,\fi{}  \beta,\ifmmode\,\fi{}  \gamma)$.
\end{dem}

\begin{rem}
  Le même résultat s'étend à des \emph{pro}groupoïdes : on peut écrire
  $A = \proobj_{\lambda \in \Lambda} A_\lambda$,\ifmmode\,\fi{} 
  $B = \proobj_{\lambda \in \Lambda} B_\lambda$ et
  $C = \proobj_{\lambda \in \Lambda} C_\lambda$
  et reprendre la même démonstration,\ifmmode\,\fi{}  en rajoutant 
  « pour $\lambda \geq \lambda_0$ »,\ifmmode\,\fi{}  
  « pour $\lambda \geq \lambda_1$ » ou
  « pour $\lambda \geq \sup\{\,\ifmmode\,\fi{}  \lambda_0,\ifmmode\,\fi{}  \lambda_1 \,\ifmmode\,\fi{} \}$ »
  un peu partout. 
\end{rem}

\begin{coro}\label{coro:VK1cocart}%
  Soient $\T$ un topos localement connexe,\ifmmode\,\fi{}  
  $U$ et $V$ des objets tels que 
  $U \amalg V \lra *$ soit un épimorphisme et tels que 
  $U \times U \iso U$ et $V \times V \iso V$. 
  Soit $P = (p_i : \cat{Ens} \lra \SLC\T)_{i\in I}$ un ensemble de
  points-bases de $U \times V$,\ifmmode\,\fi{}  tel qu'il y ait au moins un point dans
  chaque composante connexe de $U \times V$,\ifmmode\,\fi{}  $U$,\ifmmode\,\fi{}  $V$ et $*$.
  On a alors un carré cocartésien de groupoïdes fondamentaux
  $$ \xymatrix{ \relax
    \pi_1(\T/ U \times V,\ifmmode\,\fi{} P) \ar[r] \ar[d] & 
    \pi_1(\T/ U,\ifmmode\,\fi{} P) \ar[d] \\ \pi_1(\T/ V,\ifmmode\,\fi{} P) \ar[r] & 
    \pi_1(\T\!,\ifmmode\,\fi{} P). \cocartesien } $$
\end{coro}


\section{Théorème de Van Kampen pour les champs algébriques}

Nous allons maintenant appliquer le corollaire \ref{coro:VK1cocart} au
topos $\cat{Sh}\,\ifmmode\,\fi{} X$ des faisceaux étales sur un champ algébrique $X$. 

\begin{theo}\label{theo:VKchamps}%
  Soient $X$ un champ algébrique,\ifmmode\,\fi{}  $U$ et $V$ des sous-champs ouverts tels
  que $X = U \union V$ et 
  $P = (p_i : \cat{Ens} \lra \SLCF \cat{Sh}\,\ifmmode\,\fi{}  X)_{i\in I}$ un ensemble
  de points-bases de $U \dtimes_X V$,\ifmmode\,\fi{}  avec au moins un
  point dans chaque composante connexe de $U \dtimes_X V$,\ifmmode\,\fi{}  $U$,\ifmmode\,\fi{}  $V$ et
  $X$. On a alors un carré cocartésien de progroupoïdes,\ifmmode\,\fi{}  
  $$ \xymatrix{\relax
    \pi_1(U \inter V,\ifmmode\,\fi{}  P) \ar[r] \ar[d] & \pi_1(U,\ifmmode\,\fi{} P) \ar[d] \\
    \pi_1(V,\ifmmode\,\fi{} P) \ar[r] & \pi_1(X,\ifmmode\,\fi{} P). \cocartesien} $$
\end{theo}

\begin{dem}
  Comme $U \lra X$ est un sous-champ ouvert,\ifmmode\,\fi{}  
  c'est un objet de $X'\et$,\ifmmode\,\fi{}  auquel correspond un préfaisceau
  représentable $\wt U$,\ifmmode\,\fi{}  objet de 
  $\cat{Sh}\,\ifmmode\,\fi{} X'\et$.

  Montrons que 
  $\wt U \times \wt U \iso \wt U$ dans $\cat{Sh}\,\ifmmode\,\fi{} X$.
  D'après le lemme \ref{lemme:souschamp},\ifmmode\,\fi{}  on a 
  \begin{align*}
    U \Dtimes_X U &\equiv U
    &&\hspace{-3cm}\text{dans } \dcat{Champs} \\
    \intertext{donc,\ifmmode\,\fi{}  d'après le lemme \ref{lemme:2prod2prodfib},\ifmmode\,\fi{} }
    U \Dtimes U &\equiv U
    &&\hspace{-3cm}\text{dans } \dcat{Champs}/X \\
    \intertext{donc,\ifmmode\,\fi{}  d'après le lemme \ref{lemme:2prodprod},\ifmmode\,\fi{} }
    U \times U &\iso U 
    &&\hspace{-3cm}\text{dans } \cat{Cat}( \dcat{Champs}/X ). \\
    \intertext{Or,\ifmmode\,\fi{}  $X'\et$ est une sous-catégorie pleine de 
      $\cat{Cat}( \dcat{Champs}/X )$,\ifmmode\,\fi{}  donc}
    U \times U &\iso U 
    &&\hspace{-3cm}\text{dans } X'\et \\
    \intertext{donc}
    \wt{U \times U} &\iso \wt U 
    &&\hspace{-3cm}\text{dans } \cat{Sh}\,\ifmmode\,\fi{}  X'\et \\
    \intertext{donc (puisque le plongement d'Yoneda est exact à
      gauche),\ifmmode\,\fi{} }
    \wt U \times \wt U &\iso \wt U 
    &&\hspace{-3cm}\text{dans } \cat{Sh}\,\ifmmode\,\fi{}  X'\et.
  \end{align*}    

  On peut maintenant appliquer le corollaire
  \ref{coro:VK1cocart} au topos 
  $\cat{Sh}\,\ifmmode\,\fi{} X'\et$ et aux objets 
  $\wt U$ et $\wt V$ qui vérifient bien 
  $ \wt U \times \wt U \iso \wt U$ et 
  $\wt V \times \wt V \iso \wt V$ 
  et $\wt U \amalg \wt V \lra *$ est un épimorphisme (car la famille 
  $(U \lra *,\ifmmode\,\fi{}  V \lra *)$ est couvrante,\ifmmode\,\fi{}  donc 
  $(\wt U \lra *,\ifmmode\,\fi{}  \wt V \lra *)$ aussi).
  On a donc un carré cocartésien de progroupoïdes 
  $$ \xymatrix{
    \pi_1( \cat{Sh}\,\ifmmode\,\fi{}  X'\et / \wt U \times \wt V,\ifmmode\,\fi{}  P ) 
    \ar[r] \ar[d] &
    \pi_1( \cat{Sh}\,\ifmmode\,\fi{}  X'\et / \wt U ,\ifmmode\,\fi{}  P ) 
    \ar[d] \\
    \pi_1( \cat{Sh}\,\ifmmode\,\fi{}  X'\et / \wt V,\ifmmode\,\fi{}  P ) 
    \ar[r] &
    \pi_1( \cat{Sh}\,\ifmmode\,\fi{}  X'\et,\ifmmode\,\fi{}  P )
    \cocartesien
    }
  $$
  Les lemmes \ref{lemme:XetUequivUet}
  et \ref{lemme:Xet/U}  montrent que 
  $\cat{Sh}\,\ifmmode\,\fi{} X'\et/\wt U \equiv \cat{Sh}\,\ifmmode\,\fi{} U'\et$ et
  $\cat{Sh}\,\ifmmode\,\fi{} X'\et/\wt V \equiv \cat{Sh}\,\ifmmode\,\fi{} V'\et$ ; d'autre part,\ifmmode\,\fi{}  
  d'après les lemmes 
  \ref{lemme:2prod2prodfib} et \ref{lemme:2prodprod}
  et comme
  le plongement d'Yoneda est exact à gauche,\ifmmode\,\fi{}  
  $$\wt U \times \wt V \equiv \wt{ U \Dtimes_X V },\ifmmode\,\fi{}  $$
  donc 
  $$\cat{Sh}\,\ifmmode\,\fi{} X'\et/\wt U \times \wt V \equiv 
  \cat{Sh}\,\ifmmode\,\fi{} ( U \Dtimes_X V )'\et.$$
  Finalement,\ifmmode\,\fi{}  on a un carré cocartésien de progroupoïdes 
  $$ \raisebox{\depth}{\xymatrix{
      \pi_1( \cat{Sh}\,\ifmmode\,\fi{}  ( U \Dtimes_X V )'\et ,\ifmmode\,\fi{}  P ) 
      \ar[r] \ar[d] &
      \pi_1( \cat{Sh}\,\ifmmode\,\fi{}  U'\et ,\ifmmode\,\fi{}  P ) 
      \ar[d] \\
      \pi_1( \cat{Sh}\,\ifmmode\,\fi{}  V'\et,\ifmmode\,\fi{}  P ) 
      \ar[r] &
      \pi_1( \cat{Sh}\,\ifmmode\,\fi{}  X'\et,\ifmmode\,\fi{}  P ).
      \cocartesien
      }}
  \qed
  $$
\end{dem}

\begin{lemme}\label{lemme:souschamp}%
  Si $i : U \lrh X$ est un sous-champ d'un champ algébrique,\ifmmode\,\fi{}  
  alors 
  $U \equiv U \Dtimes_X U$.
\end{lemme}
\begin{dem}
  Le sous-champ $i : U \lrh X$ est une sous-catégorie
  fibrée. La fibre au dessus d'un schéma $T$ de la catégorie fibrée 
  $U \dtimes_X U$ a pour objets les $(f,\ifmmode\,\fi{} g,\ifmmode\,\fi{} \alpha)$,\ifmmode\,\fi{}  où 
  $f \in U(T)$,\ifmmode\,\fi{}  $g \in U(T)$ et 
  $\alpha \in \hom_{X(T)}( if,\ifmmode\,\fi{}  ig) = 
  \hom_{U(T)}(f,\ifmmode\,\fi{}  g)$,\ifmmode\,\fi{}  
  \ie,\ifmmode\,\fi{}  les diagrammes 
  $$ \xymatrix{ T \UN fg\alpha & U,\ifmmode\,\fi{}  } $$
  et les morphismes
  $(f,\ifmmode\,\fi{} g,\ifmmode\,\fi{} \alpha) \lra (f',\ifmmode\,\fi{} g',\ifmmode\,\fi{} \alpha')$ 
  sont les $(\phi,\ifmmode\,\fi{} \psi)$ où 
  $$ \xymatrix{
    f \ar[d]_{\phi} & if \ar[d]_{i\phi} \ar[r]^\alpha & ig \ar[d]^{i\psi} & g
    \ar[d]^\psi \\
    f' & if' \ar[r]_{\alpha'} & ig' & g'. } $$
  Mais tout $(f,\ifmmode\,\fi{} g,\ifmmode\,\fi{} \alpha)$ est isomorphe à 
  $(f,\ifmmode\,\fi{}  \id_f,\ifmmode\,\fi{}  f)$ : 
  $$ \xymatrix{
    f \ar[d]_{\id_f} & f \ar[d]_{\id_f} \ar[r]^\alpha & g
    \ar[d]^{\alpha^{-1}} & g \ar[d]^{\alpha^{-1}} \\
    f & f \ar[r]_{\id_f} & f & f. } $$
  D'autre part les morphismes 
  $(f,\ifmmode\,\fi{} f,\ifmmode\,\fi{} \id) \lra (g,\ifmmode\,\fi{} g,\ifmmode\,\fi{} \id)$ sont les $(\phi,\ifmmode\,\fi{} \psi)$ tels que 
  $$ \xymatrix{ 
    f \ar[d]_\phi & f \ar[d]_\phi \ar[r]^{\id} & f \ar[d] ^\psi & f
    \ar[d] ^\psi \\
    g & g \ar[r]_{\id_g} & g & g,\ifmmode\,\fi{}  } $$
  \ie,\ifmmode\,\fi{}  tels que $\phi = \psi$. 
  Cela donne une équivalence de catégories 
  $$ U_T \lra \left( 
    U \Dtimes_X U \right)_T. $$
\end{dem}

\begin{lemme}\label{lemme:XetUequivUet}%
  Si $U \lrh X$ est un sous-champ,\ifmmode\,\fi{}  
  on a une équivalence de sites
  $X'\et / U \equiv U'\et$.
\end{lemme}
\begin{dem}
  (a) Les objets de 
  $U'\et$ sont les 
  $t : T \lra U$ et les morphismes 
  $(T,\ifmmode\,\fi{} t) \lra (T',\ifmmode\,\fi{}  t')$ sont les 
  \cbstart
  $(f:T \lra T',\ifmmode\,\fi{}  \alpha:t'f\lra t)$,\ifmmode\,\fi{} 
  \cbend
  $$\english \xymatrix{ 
    T \ar[rr] ^f \ar[rd] |*{}="A" _t 
    \enlargexyentry A
    && T' \ar[ld]^{t'} 
    \ar@{=>}"A"_(.7)\alpha \\ & U,} $$
  modulo les 
  $ \xymatrix{ T \UN{f}{f'}{\theta} & T' } $
  tels que 
  $\raisebox{.5\depth}{\xymatrix{ t'f \ar[r]^{t' \theta} \ar[rd]_\alpha & t'f'
    \ar[d]^{\alpha'} \\ & t. }} $

  (b) Les objets de 
  $X'\et/U$ sont les 
  \cbstart
  $(T,\ifmmode\,\fi{} s:T\lra X,\ifmmode\,\fi{} t:T\lra U,\ifmmode\,\fi{} \beta:ut\lra s)$,\ifmmode\,\fi{} 
  \cbend
  $$\english \xymatrix{ 
    T \ar[rr] ^t \ar[rd] |*{}="A" _s 
    \enlargexyentry A
    && U \ar[ld]^{u} 
    \ar@{=>}"A"_(.7)\beta \\ & X,} $$
  modulo les 
  $\xymatrix{ T \UN{t}{t'}{\phi} & U } $
  tels que 
  $\raisebox{.5\depth}{\xymatrix{ s't \ar[r]^{s'\phi} \ar[rd]_\beta & s't'
    \ar[d]^{\beta'} \\ & s}} $
  et les morphismes 
  $(T,\ifmmode\,\fi{} s,\ifmmode\,\fi{} t,\ifmmode\,\fi{} \beta) \lra (T',\ifmmode\,\fi{} s',\ifmmode\,\fi{} t',\ifmmode\,\fi{} \beta')$ sont les 
  \cbstart
  $(f:T\lra T',\ifmmode\,\fi{}  \alpha:t'f\lra t,\ifmmode\,\fi{}  \gamma:s'f\lra
  s)$
  tels que
  $$ \english
  \raisebox{.5\depth}{%
    \xymatrix{ 
      T \ar[rr]^f 
      \ar[rd] |*{}="A" _s 
      \enlargexyentry A
      && T' \ar[ld]^{s'} 
      \ar@{=>}"A"_(.8)\gamma \\
      & X}}
  =
  \raisebox{.5\depth}{\xymatrix{ 
      T \ar[rr]^f \ar[rd] |*{}="A" _t 
      \ar@/_1pc/[ddr] |*{}="B" _s 
      \enlargexyentry A
      && T' \ar[ld] ^{t'}
      \ar@{=>}"A"_(.8)\alpha
      \ar@/^1pc/[ldd] |*{}="C" ^{s'} 
      \enlargexyentry B
      \enlargexyentry C
      \\ & U \ar[d]_u 
      \ar@{=>}"B"^\beta 
      \ar@{=>}"C"_{\beta'}
      \\ & X,\ifmmode\,\fi{} 
      }}
  \qquad
  \text{\ie,\ifmmode\,\fi{} }
  \qquad
  \gamma = \beta \alpha \beta^{-1}
  $$
  \cbend
  modulo les 
  $ \xymatrix{ T \UN{f}{f'}{\theta} & T' } $
  tels que 
  $\raisebox{.5\depth}{\xymatrix{ s'f \ar[r]^{s'\theta} \ar[rd]_\beta & s'f
    \ar[d]^{\beta'} \\ & s }} $
  et modulo un choix différent de $\alpha$. 

  (c) Nous allons montrer que le foncteur suivant est une équivalence
  de catégories. 
  $$\english 
  \definefunctor{U'\et}{X'\et/U}%
  {\raisebox{.5\depth}{\xymatrix{T \ar[d]_t \\ U}}}%
  {\raisebox{.5\depth}{\xymatrix{ T \ar[rr]^t |*{}="A"\ar[rd]_{ut} 
        \enlargexyentry A 
        &&
        U \ar[ld]^u \\ & X \ar@{}"A"|{\circlearrowleft} 
        }}}%
  {\raisebox{.5\depth}{\xymatrix{ T \ar[rr]^f \ar[rd] |*{}="A" _t 
        \enlargexyentry A
        && T' \ar[ld]^{t'}
        \ar@{=>}"A"_(.7)\alpha \\ & U }}}%
  {\raisebox{.5\depth}{%
    \xymatrix{ 
      T \ar[rr]^f \ar[rd] |*{}="A" _t 
      \ar@/_/[ddr] |*{}="B" _{ut}
      \enlargexyentry A
      \enlargexyentry B 
      && T' \ar[ld] ^{t'} \ar@{=>}"A"_(.7)\alpha
      \ar@/^/[ldd] |*{}="C" ^{ut'}
      \\ & U \ar[d]_u 
      \ar@{}"B"^{\circlearrowleft}
      \ar@{}"C"_{\circlearrowleft}
      \\ & X 
      }}}%
  $$
  
  (d) Il est essentiellement surjectif,\ifmmode\,\fi{}  car tout objet 
  $(T,\ifmmode\,\fi{} s,\ifmmode\,\fi{} t,\ifmmode\,\fi{} \beta)$ est isomorphe à un objet de la forme 
  $(T,\ifmmode\,\fi{}  ut,\ifmmode\,\fi{}  t,\ifmmode\,\fi{}  \id)$ 
  par l'isomorphisme 
  $$ \xymatrix{ 
    (T,\ifmmode\,\fi{}  s,\ifmmode\,\fi{}  t,\ifmmode\,\fi{}  \beta) \ar[rr]^{ (\id,\ifmmode\,\fi{}  \id,\ifmmode\,\fi{}  \beta) } &&
    (T,\ifmmode\,\fi{}  ut,\ifmmode\,\fi{}  t,\ifmmode\,\fi{}  \id). }$$
  
  (e) Montrons maintenant qu'il est pleinement fidèle. 
  Soient donc 
  $(T,\ifmmode\,\fi{}  t)$ et $(T',\ifmmode\,\fi{}  t')$ des objets de $U'\et$. Montrons que
  l'application 
  $$ \definefunction
  {\hom \bigl( 
    (T,\ifmmode\,\fi{}  t),\ifmmode\,\fi{}  
    (T',\ifmmode\,\fi{}  t') 
    \bigr)}%
  {\hom\bigl( 
    (T,\ifmmode\,\fi{}  ut,\ifmmode\,\fi{}  u,\ifmmode\,\fi{}  \id),\ifmmode\,\fi{}  
    (T,\ifmmode\,\fi{}  ut',\ifmmode\,\fi{}  u,\ifmmode\,\fi{}  \id)
    \bigr)}%
  {(f,\ifmmode\,\fi{} \alpha)}{(f,\ifmmode\,\fi{} \alpha,\ifmmode\,\fi{} u\alpha)}
  $$ est bijective. La surjectivité est imédiate,\ifmmode\,\fi{}  car tout morphisme 
  $(f,\ifmmode\,\fi{}  \alpha,\ifmmode\,\fi{}  \gamma)$ est l'image du morphisme 
  $(f,\ifmmode\,\fi{}  \alpha)$,\ifmmode\,\fi{}  car $\gamma$ s'exprime en fonction de $f$ et
  $\alpha$. Pour l'injectivité,\ifmmode\,\fi{}  il suffit de remarquer que 
  $\alpha$ est entièrement déterminé par la donnée de $f$ et $u
  \alpha$ : si 
  $$ 
  \english
  \raisebox{.5\depth}{%
    \xymatrix{ T \ar[rr]^f \ar[rd] |*{}="A" 
      \enlargexyentry A
      && T' \ar[ld] 
      \ar@{=>}"A"_(.7){\alpha} 
      \\ & U
      \ar[d]^u \\ &X }}
  =
  \raisebox{.5\depth}{%
    \xymatrix{ T \ar[rr]^f \ar[rd] |*{}="A" 
      \enlargexyentry A
      && T' \ar[ld] 
      \ar@{=>}"A"_(.7){\alpha'} 
      \\ & U
      \ar[d]^u \\ &X,\ifmmode\,\fi{}  }}
  $$
  alors $\alpha = \alpha'$,\ifmmode\,\fi{}  car le foncteur 
  $\cat{hom}(T,\ifmmode\,\fi{} U) \lra \cat{hom}(T,\ifmmode\,\fi{} X)$ est pleinement fidèle,\ifmmode\,\fi{}  
  car~$U$ est un sous-champ.
  
  (f) Enfin,\ifmmode\,\fi{}  il s'agit d'une équivalence de sites,\ifmmode\,\fi{}  car les topologies
  sont définies de la même manière : les familles étales surjectives
  (ce sont les mêmes dans $X'\et/U$ et $U'\et$) sont décrétées
  couvrantes.
\end{dem}

\begin{lemme}\label{lemme:Xet/U}%
  Soit $U \lra X$ un sous-champ  ouvert
  et $\wt U$ le faisceau étale sur~$X$ correspondant. 
  On a une équivalence de catégories 
  $\cat{Sh}\,\ifmmode\,\fi{}  X\et / \wt U \equiv 
  \cat{Sh}\,\ifmmode\,\fi{}  U\et.$
\end{lemme}
\begin{dem}
  \begin{align*}
    \cat{Sh}\,\ifmmode\,\fi{}  X\et / \wt U 
    &\equiv 
    \cat{Sh}\,\ifmmode\,\fi{}  X'\et / \wt U  
    &&\hspace{-2cm}\text{d'après le lemme \ref{lemme:deuxsites}} 
    \\
    &\equiv 
    \cat{Sh} ( X'\et / U )
    &&\hspace{-2cm}\text{d'après \cite[III.5.4]{SGA4}}
    \\
    &\equiv
    \cat{Sh}\,\ifmmode\,\fi{}  U'\et 
    &&\hspace{-2cm}\text{d'après le lemme \ref{lemme:XetUequivUet}}
    \\
    &\equiv 
    \cat{Sh}\,\ifmmode\,\fi{}  U\et
    &&\hspace{-2cm}\text{d'après le lemme \ref{lemme:deuxsites}.}
  \end{align*}
\end{dem}

\cbstart
\bibliographystyle{hubris}%
\nocite{*}%
\bibliography{1}
\cbend

\end{document}